\documentclass[11pt]{article}

\usepackage{amsmath, amsfonts}
\usepackage{adjustbox}
\usepackage{microtype}
\usepackage{graphicx}
\usepackage{arydshln}
\usepackage{fancyvrb}
\usepackage{hyperref}
\usepackage{fullpage} 
\usepackage{caption}
\usepackage{cite}

\title{Constructing Optimal Kobon Triangle Arrangements via Table Encoding, SAT Solving, and Heuristic Straightening}
\author{
  Pavlo Savchuk \\
  \texttt{zegalur@gmail.com}
}
\date{July 10, 2025} 

\begin{document}

\maketitle

\begin{abstract}
We present new methods and results for constructing optimal Kobon triangle arrangements. First, we introduce a compact table notation for describing arrangements of pseudolines, enabling the representation and analysis of complex cases, including symmetrical arrangements, arrangements with parallel lines, and arrangements with multiple-line intersection points. Building on this, we provide a simple heuristic method and tools for recovering straight-line arrangements from a given table, with the ability to enforce additional properties such as symmetries. The tool successfully recovers arrangements for many previously known optimal solutions. Additionally, we develop a tool that transforms the search for optimal Kobon arrangement tables into a SAT problem, allowing us to leverage modern SAT solvers (specifically Kissat) to efficiently find new solutions or to show that no other solutions exist (for example, confirming that no simple perfect arrangement with 33 triangles exists in the 11-line case). Using these techniques, we find new optimal Kobon arrangements for 23 and 27 lines, along with several other new results.
\end{abstract}

\section{Introduction}
\textls[-10] {
The classical Kobon triangle problem asks for the largest number $N(n)$ of nonoverlapping triangles that can be constructed using $n$ straight lines on a plane \cite{wiki:kobon,wolfram:kobon}. As the problem remains unsolved, tight upper bounds on the values of $N(n)$ are known \cite{cb2007a, bartholdi2007}. One way to find $N(n)$ for a specific $n$ is to construct an (optimal) arrangement that meets a known upper bound (see~\cite[A006066]{oeis}).

Instead of constructing arrangements with straight lines, a common strategy is to loosen the constraints and look for optimal arrangements of pseudolines instead. This task is easier to formalize, because we do not care about linearity, but the downside is that after finding a candidate arrangement, we then need to ``straighten'' it by constructing a corresponding arrangement of straight lines. 

Another common strategy is to utilize symmetries, specifically mirror and rotational symmetries. It is often a good idea to first search for a symmetrical solution because it is combinatorially simpler, and many examples of symmetrical optimal arrangements are already known for different $n$.

In this paper, we show how to automate these strategies by utilizing a table notation, a SAT solver, and heuristic straightening. The results are provided as ready-to-use tools available at \cite{koboncnf, lineorder}. Some of the results are also included at the end of this paper in the Appendix~\ref{appendix:tan_form}-\ref{appendix:big_ar}.

There are other approaches, equivalent notions, and relevant mathematical structures that aren’t discussed here (such as \textit{abstract order types}\cite[Section~3 and Appendix~B]{rote2025numpslaexperimentalresearch} or \textit{oriented matroids}\cite[Section~2]{Fukuda2013}, etc.), as well as advanced tools for assisting with the research\cite{rote2025numpslaexperimentalresearch, Rote_NumPSLA_2024}. For more information, please refer to the list of sources at the end of this paper and their references. 
}

\section{Table Notation}

This section introduces a simple method for encoding complex arrangements of lines as tables. The idea is to develop a notation that is expressive enough to represent optimal Kobon arrangements (including those with parallel lines and multiple-line intersection points) in a more or less invariant form, similar to \textit{wiring diagrams}\cite{geom2017wiring}.

\subsection{Construction}

\begin{figure}
\centering
\begin{tabular}{p{5cm}p{5cm}p{5cm}}
\begin{minipage}[t]{\linewidth}
\centering
\includegraphics[width=\textwidth]{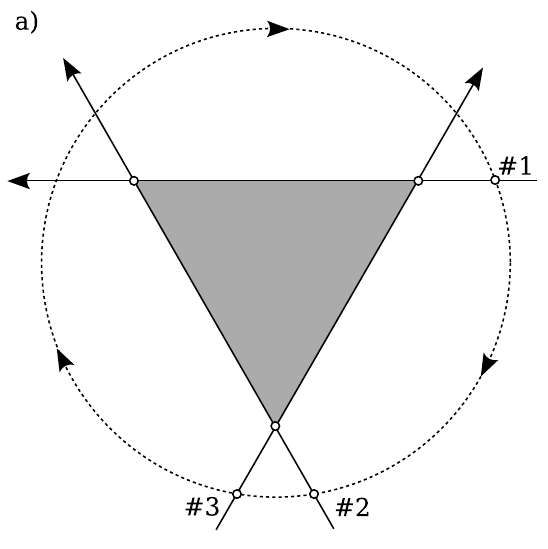} \\
\begin{BVerbatim}
[[3,2],
 [3,1],
 [2,1]]
\end{BVerbatim}
\end{minipage}
&
\begin{minipage}[t]{\linewidth}
\centering
\includegraphics[width=\textwidth]{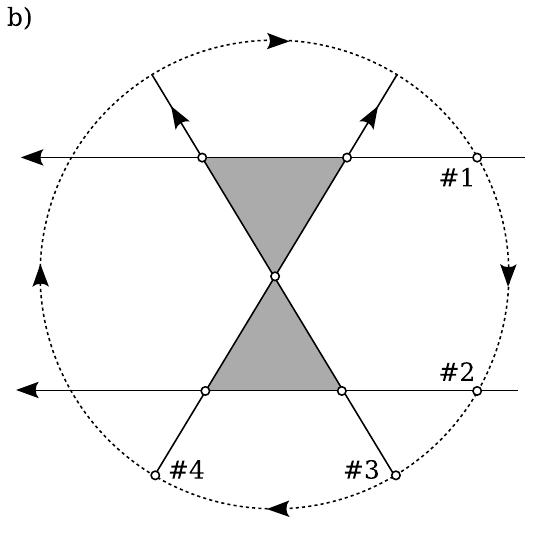} \\
\begin{BVerbatim}
[[4,3],
 [3,4],
 [2,4,1],
 [2,3,1]]
\end{BVerbatim}
\end{minipage}
&
\begin{minipage}[t]{\linewidth}
\centering
\includegraphics[width=\textwidth]{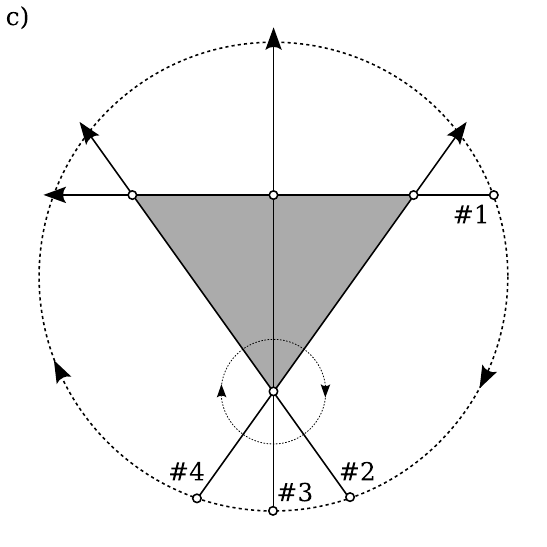} \\
\begin{BVerbatim}
[[4,3,2],
 [[3,4],1],
 [[4,2],1],
 [[2,3],1]]
\end{BVerbatim}
\end{minipage}
\end{tabular}

\caption{Some simpler arrangements and their tables. In this paper, we always use the topmost horizontal line as line~\#1, with its direction from right to left. Note that the $i$-th row corresponds to the $i$-th line. (a)~Three lines forming a triangle. (b)~Four lines, two of them parallel. (c)~An arrangement with a multiple-line intersection point.}
\label{fig:tab_examples}
\end{figure}

To build a table from an arrangement, we first enclose all the crossing points of the arrangement within a large circle. We then select one line to be line~\#1. We choose the point where line~\#1 enters the circle (typically, the rightmost point). From this point, we proceed clockwise from line~\#1 around the circle and number each other line that enters the circle as line~\#2, line~\#3, and so on. Already numbered lines are skipped if needed. We continue this process until all the lines are numbered. In this way, all lines are enumerated, and each line receives an assigned direction (e.g, see Figure~\ref{fig:tab_examples}). For arrangements with parallel lines, we select line~\#1 so that groups of parallel lines are given consecutive numbers (see part {(b)} of Figure~\ref{fig:tab_examples}).

Finally, for every line, from \#1 to \#$n$, we write down a row of line numbers representing the lines that cross it, in the order determined by its assigned direction. For intersections involving more than two lines, we write a group of all intersecting lines in clockwise order, starting from the current line. In this way, we construct an arrangement table. Some simpler arrangements with their numbered lines and corresponding tables are shown in Figure~\ref{fig:tab_examples}. In this paper, we always use the topmost horizontal line as line \#1, with its direction from right to left. Note that the $i$-th row corresponds to the $i$-th line.

\subsection{Properties}

\begin{figure}
\centering
\includegraphics[width=0.8\textwidth]{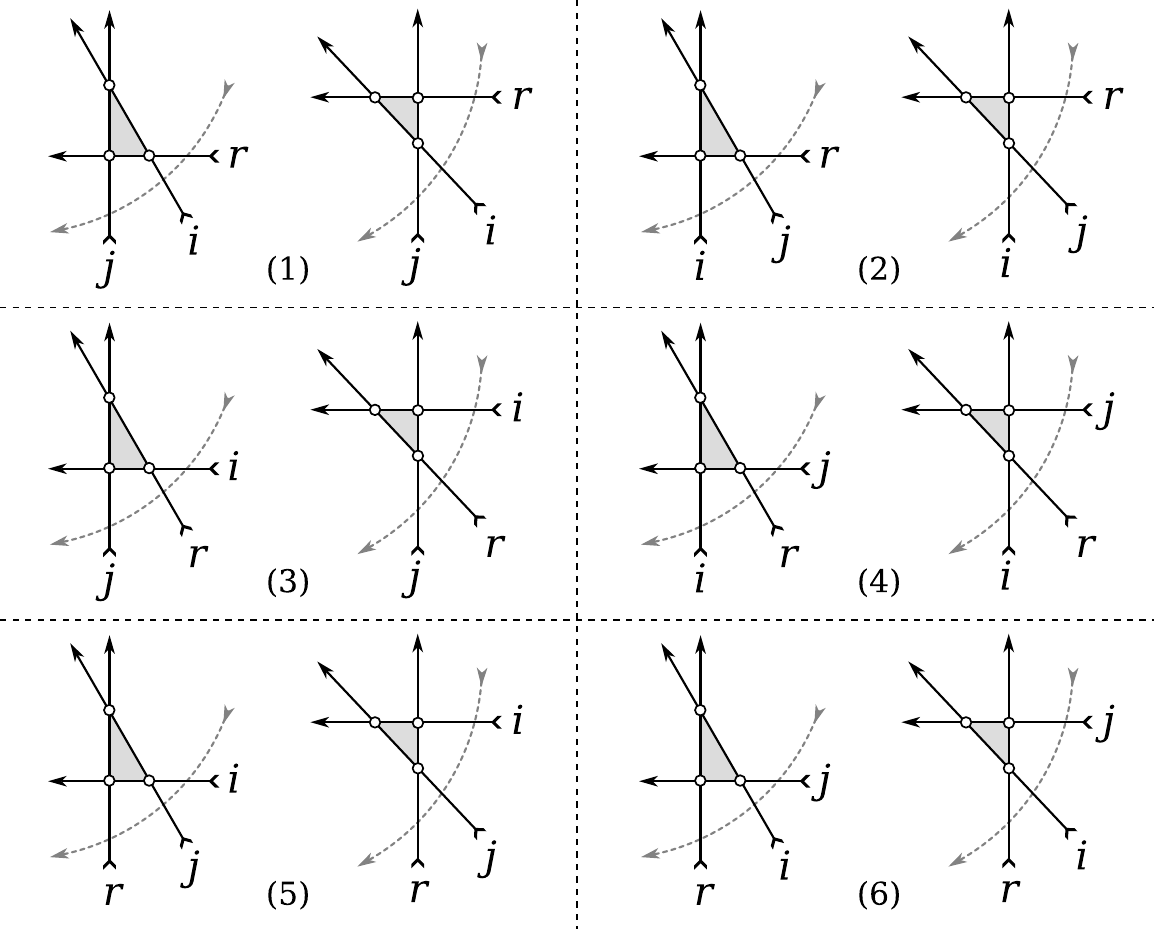}
\caption{For three different lines $r$, $i$, and $j$, there are only twelve possible ways they can form a triangle (grouped into six classes by the ordering: $rij$, $rji$, $\dots$, $jir$).}
\label{fig:consistency}
\end{figure}

Because tables are restricted by the logic of how (pseudo)lines can behave, they have many useful properties. Only some of these properties are listed here; an exhaustive and formalized set of properties for optimal Kobon arrangements with $n \bmod 6 \in \{3,5\}$ lines is provided in the CNF-model section. Some general properties are:
\begin{enumerate}
\item \label{prop:selfexclusion} \textit{(Self-exclusion)} 
Line~$k$ cannot appear in row~$k$.
\item \label{prop:nonrep} \textit{(Non-repetition)} 
Each line can appear at most once in a given row.
\item \label{prop:complitness} \textit{(Completeness)} 
For arrangements where all lines pairwise intersect, each row~$k$ contains all the lines $1,2,\dots,n$ except $k$ itself.
\item \label{prop:consistency} \textit{(Consistency)} 
The relative ordering of any two lines $i$ and $j$ (excluding multiple-line intersection points) in row~$r$ is consistent with how (pseudo)line intersections behave in the triangle formed by lines $i$, $j$, and $r$ (see Figure~\ref{fig:consistency}). For example, if $r < i < j$ and $i$ appears before $j$ in row~$r$, then $r$ appears before $j$ in row~$i$, and $r$ appears before $i$ in row~$j$.
\end{enumerate}
An arrangement is \textit{simple} if no two lines are parallel and no three lines meet at a point. It is \textit{perfect} if every finite segment is a side of exactly one triangle and it meets the upper bound (1) from~\cite{cb2007a}.

\section{CNF Model}

\begin{figure}
\centering
\begin{tabular}{p{5cm}p{5cm}p{5cm}}
\begin{minipage}[t]{\linewidth}
\centering
\includegraphics[width=\textwidth]{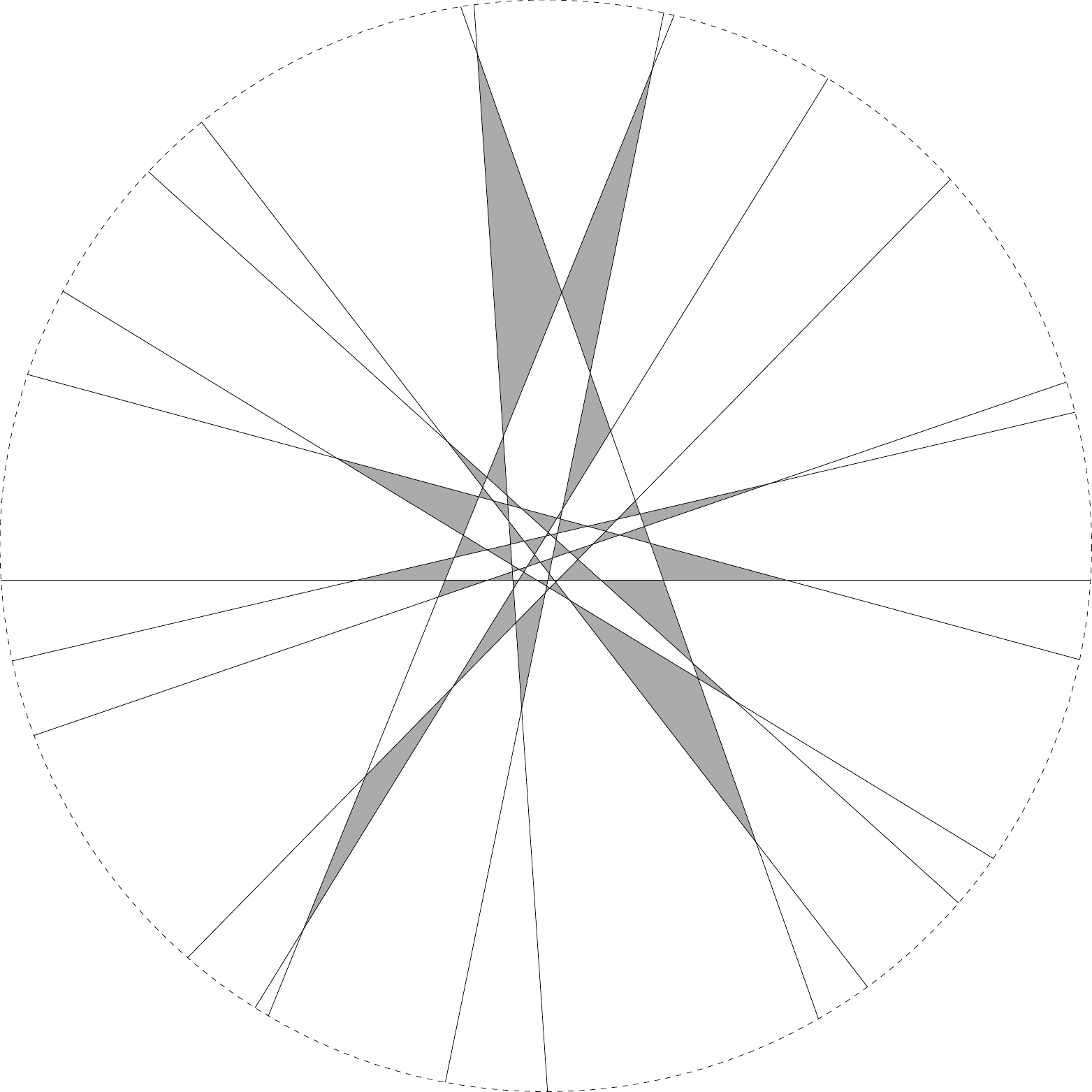} \\
\tiny
\begin{BVerbatim}

[[2,6,4,11,5,8,3,10,7,12,9,13],
 [1,6,12,11,13,8,10,4,7,5,9,3],
 [4,6,5,11,8,1,10,12,7,13,9,2],
 [3,6,1,11,12,8,13,10,2,7,9,5],
 [6,3,11,1,8,12,10,13,7,2,9,4],
 [5,3,4,1,2,12,13,11,10,8,9,7],
 [8,11,10,1,12,3,13,5,2,4,9,6],
 [7,11,3,1,5,12,4,13,2,10,6,9],
 [10,11,12,1,13,3,2,5,4,7,6,8],
 [9,11,7,1,3,12,5,13,4,2,8,6],
 [9,10,7,8,3,5,1,4,12,2,13,6],
 [9,1,7,3,10,5,8,4,11,2,6,13],
 [1,9,3,7,5,10,4,8,2,11,6,12]]
\end{BVerbatim}
\normalsize
\end{minipage}
&
\begin{minipage}[t]{\linewidth}
\centering
\includegraphics[width=\textwidth]{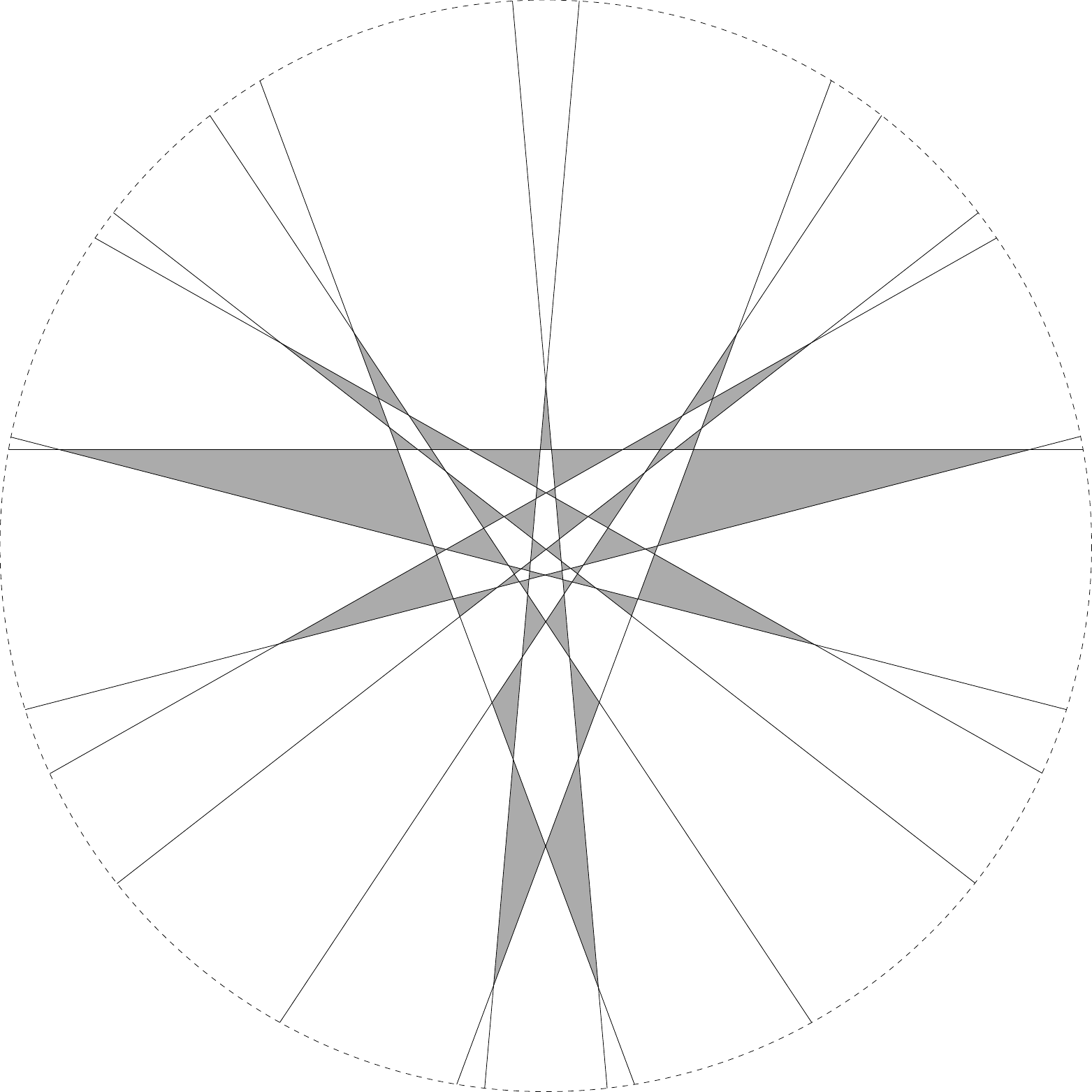} \\
\tiny
\begin{BVerbatim}

[[13,9,11,10,12,7,8,3,5,4,6,2],
 [3,9,4,10,7,13,8,11,5,12,6,1],
 [2,9,13,10,11,7,12,8,1,5,6,4],
 [9,2,10,13,7,11,8,12,5,1,6,3],
 [9,7,10,8,13,11,2,12,4,1,3,6],
 [7,9,8,10,11,13,12,2,1,4,3,5],
 [6,9,5,10,2,13,4,11,3,12,1,8],
 [9,6,10,5,13,2,11,4,12,3,1,7],
 [8,6,7,5,4,2,3,13,1,11,12,10],
 [6,8,5,7,2,4,13,3,11,1,12,9],
 [6,13,5,2,8,4,7,3,10,1,9,12],
 [13,6,2,5,4,8,3,7,1,10,9,11],
 [12,6,11,5,8,2,7,4,10,3,9,1]]
\end{BVerbatim}
\normalsize
\end{minipage}
&
\begin{minipage}[t]{\linewidth}
\centering
\includegraphics[width=\textwidth]{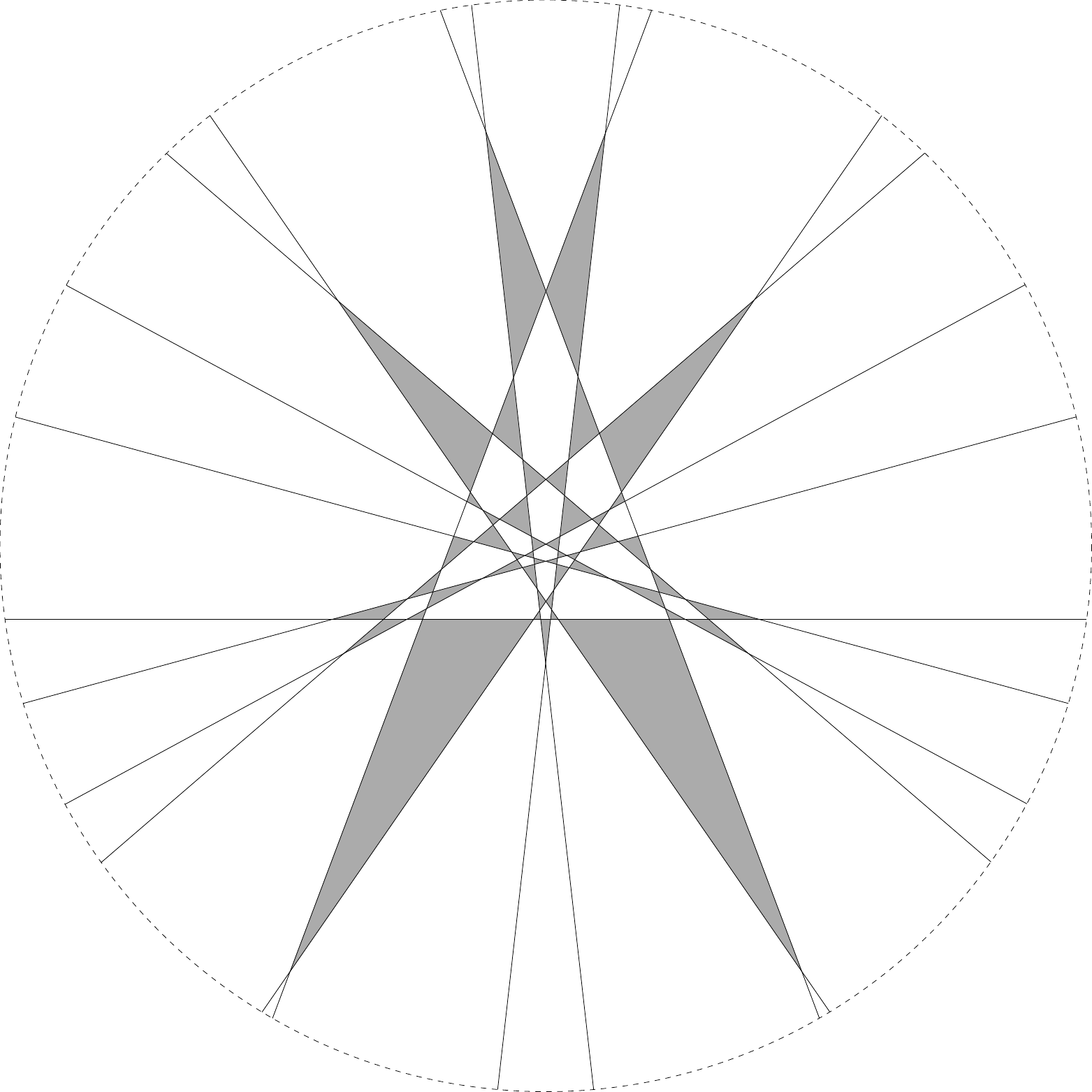} \\
\tiny
\begin{BVerbatim}

[[2,4,3,6,5,8,7,10,9,12,11,13],
 [1,4,6,3,10,8,13,7,12,5,11,9],
 [4,1,6,2,10,13,8,12,7,11,5,9],
 [3,1,2,6,13,10,12,8,11,7,9,5],
 [6,1,8,10,7,13,12,2,11,3,9,4],
 [5,1,3,2,4,13,12,10,11,8,9,7],
 [8,1,10,5,13,2,12,3,11,4,9,6],
 [7,1,5,10,2,13,3,12,4,11,6,9],
 [10,1,12,13,11,2,3,5,4,7,6,8],
 [9,1,7,5,8,2,3,13,4,12,6,11],
 [12,1,13,9,2,5,3,7,4,8,6,10],
 [11,1,9,13,5,2,7,3,8,4,10,6],
 [1,11,9,12,5,7,2,8,3,10,4,6]]
\end{BVerbatim}
\normalsize
\end{minipage}
\end{tabular}

\caption{Every unique optimal arrangement for $n=13$ where line~\#6 and line~\#9 lack a triangle, as found by \texttt{kobon-cnf}\cite{koboncnf}, shown with their respective tables.}
\label{fig:kobon-13}
\end{figure}

This section describes optimal Kobon arrangement tables for certain values of~$n$, formalized as CNF models suitable for modern SAT solvers. Specifically, when $n \bmod 6 \in \{3,5\}$, arrangements meeting the upper bound $N(n) = \frac{1}{3}n(n-2)$ are stated in~\cite{cb2007a} to be perfect, and hence simple. These properties simplify the model. We can further extend the approach to other values of~$n$ by adding segments with missing triangles (see Figure~\ref{fig:kobon-13}, etc.). A fully working CNF model implementation can be found at~\cite{koboncnf}. As SAT solvers become increasingly powerful, more results can be achieved through proper encoding (see, for example, \cite{kirchweger_et_al:LIPIcs.SAT.2023.14}).

\subsection{Variables}

In order to encode tables as a CNF statement, we define the following Boolean variables:\begin{enumerate}
\item $G(r,i,j)$ -- line $j$ is immediately after line $i$ in row~$r$ 
(for $r \neq i \land r \neq j \land i \neq j$).
\item $X(r,i,j)$ -- line $j$ is somewhere after line $j$ in row~$r$
(for $r \neq i \land r \neq j \land i \neq j$).
\item $A(r,i,k)$ -- line $i$ is in column~$k$ of row~$r$
(for $r \neq i$).
\item \textit{(Optional)} $M_k(r_k,i,j)$ -- segment $(i,j)$ on line~$r_k$ is missing a triangle.
\end{enumerate}
Here, lines are numbered from $1$ to $n$, rows from $1$ to $n$, and columns from $1$ to $n-1$. In total there are $O(n^3)$, or $\approx(3n^3-8n^2+5n)$ variables (see Table~\ref{tab:koboncnf-stats} for concrete numbers).

\subsection{Clauses}

Each rule is described textually, formalized in first-order logic, and expressed in CNF:

\begin{enumerate}

\item Each row contains all the other lines:
$$
\forall r,i \neq r: \exists k: A(r,i,k)
\ \ \ \ \ \Leftrightarrow \ \ \ \ \ 
\bigwedge\limits_{r, i \neq r} \bigvee\limits_k A(r,i,k)
$$

\item Rows do not have duplicate entries:
$$
\forall r,i \neq r,k: A(r,i,k) \rightarrow \forall j \neq k: \neg A(r,i,j)
\ \ \ \ \ \Leftrightarrow \ \ \ \ \ 
\bigwedge\limits_{r,i \neq r,k} \bigwedge\limits_{j \neq k} \neg A(r,i,k) \lor \neg A(r,i,j)
$$

\item Connection between $A$ (at) and $G$ (immediately after):
\begin{enumerate}
\item $
A(r,i,k) \land A(r,j,k+1) \rightarrow G(r,i,j)
\ \ \ \ \ \Leftrightarrow \ \ \ \ \ 
\bigwedge\limits_{\substack{r,i \neq r\\ j \neq r, j \neq i\\ k<n-1}} 
\neg A(r,i,k) \lor \neg A(r,j,k+1) \lor G(r,i,j)
$
\item $
A(r,i,k) \land G(r,i,j) \rightarrow A(r,j,k+1)
\ \ \ \ \ \Leftrightarrow \ \ \ \ \ 
\bigwedge\limits_{\substack{r,i \neq r\\ j \neq r, j \neq i\\ k<n-1}} 
\neg A(r,i,k) \lor A(r,j,k+1) \lor \neg G(r,i,j)
$
\item $
G(r,i,j) \land A(r,j,k+1) \rightarrow A(r,i,k)
\ \ \ \ \ \Leftrightarrow \ \ \ \ \ 
\bigwedge\limits_{\substack{r,i \neq r\\ j \neq r, j \neq i\\ k<n-1}} 
A(r,i,k) \lor \neg A(r,j,k+1) \lor \neg G(r,i,j)
$
\end{enumerate}

\item Only one line can immediately follow another:
\small
$$
\forall r,i \neq r, j \neq r, i \neq j: G(r,i,j) \rightarrow \forall k\neq j: \neg G(r,i,k)
\ \ \ \ \ \Leftrightarrow \ \ \ \ \ 
\bigwedge\limits_{\substack{r,i \neq r\\ j \neq r, i \neq j\\k \neq r, k \neq j, k \neq i}} 
\neg G(r,i,j) \lor \neg G(r,i,k)
$$
\normalsize

\item Connection between ``first'' and ``immediately after'':
\begin{enumerate}
\item $
A(r,i,1) \rightarrow \forall j \neq i, j \neq r: \neg G(r,j,i)
\ \ \ \ \ \Leftrightarrow \ \ \ \ \ 
\bigwedge\limits_{r,i \neq r, j \neq i, j \neq r} 
\neg A(r,i,1) \lor \neg G(r,j,i)
$
\item $
\neg A(r,i,1) \rightarrow \exists j \neq i, j \neq r: G(r,j,i)
\ \ \ \ \ \Leftrightarrow \ \ \ \ \ 
\bigwedge\limits_{r,i \neq r} 
A(r,i,1) \lor
\bigvee\limits_{\substack{j \neq i\\ j \neq r}} G(r,j,i)
$
\end{enumerate}

\item Relation between ``last'' and ``immediately after'':
\begin{enumerate}
\item $
A(r,i,n-1) \rightarrow \forall j \neq i, j \neq r: \neg G(r,i,j)
\ \ \ \ \ \Leftrightarrow \ \ \ \ \ 
\bigwedge\limits_{r,i \neq r, j \neq i, j \neq r} 
\neg A(r,i,n-1) \lor \neg G(r,i,j)
$
\item $
\neg A(r,i,n-1) \rightarrow \exists j \neq i, j \neq r: G(r,i,j)
\ \ \ \ \ \Leftrightarrow \ \ \ \ \ 
\bigwedge\limits_{r,i \neq r} 
A(r,i,n-1) \lor
\bigvee\limits_{\substack{j \neq i\\ j \neq r}}  G(r,i,j)
$
\end{enumerate}

\item $X(r,i,j)$ ($j$ somewhere after $i$) is \textit{total} (connected, complete):
$$
\forall r,i,j (i \neq r, j \neq r, i \neq j): X(r,i,j) \lor X(r,j,i)
\ \ \ \ \ \Leftrightarrow \ \ \ \ \ 
\bigwedge\limits_{\substack{r,i,j\\(i \neq r, j \neq r, i \neq j)}}
X(r,i,j) \lor X(r,j,i)
$$

\item $X(r,i,j)$ ($j$ somewhere after $i$) is \textit{antisymmetric}:
$$
\forall r,i,j (i \neq r, j \neq r, i \neq j): X(r,i,j) \rightarrow \neg X(r,j,i)
\ \ \ \ \ \Leftrightarrow \ \ \ \ \ 
\bigwedge\limits_{\substack{r,i,j\\(i \neq r, j \neq r, i \neq j)}}
\neg X(r,i,j) \lor \neg X(r,j,i)
$$

\item $X(r,i,j)$ ($j$ somewhere after $i$) is \textit{transitive}:
\begin{align*}
\forall r,i,j,k (i \neq r, j \neq r, k \neq r, i \neq j, i \neq k, j \neq k):
X(r,i,j) \land X(r,j,k) \rightarrow X(r,i,k) \ \ \ \ \ \Leftrightarrow 
\\
\Leftrightarrow \ \ \ \ \ 
\bigwedge\limits_{
\substack{r,i,j,k\ (i \neq r, j \neq r, k \neq r, i \neq j, i \neq k, j \neq k)}}
\neg X(r,i,j) \lor \neg X(r,j,k) \lor X(r,i,k)
\end{align*}

\item Relationship between $X(r,i,j)$ ($j$ somewhere after $i$) and ``first'' and ``last'' in a row:
\small
\begin{enumerate}
\item $
A(r,i,1) \rightarrow \forall j \neq i: X(r,i,j)
\ \ \ \ \ \Leftrightarrow \ \ \ \ \ 
\bigwedge\limits_{\substack{r,i,j\\(i \neq r, j \neq r, i \neq j)}}
\neg A(r,i,1) \lor X(r,i,j)
$
\item $
A(r,i,n-1) \rightarrow \forall j \neq i: X(r,j,i)
\ \ \ \ \ \Leftrightarrow \ \ \ \ \ 
\bigwedge\limits_{\substack{r,i,j\\(i \neq r, j \neq r, i \neq j)}}
\neg A(r,i,n-1) \lor X(r,j,i)
$
\item $
\bigl(\ \forall j \neq i: X(r,i,j)\ \bigr) \rightarrow A(r,i,1)
\ \ \ \ \ \Leftrightarrow \ \ \ \ \ 
\bigwedge\limits_{\substack{r,i\\(i \neq r)}}
A(r,i,1) \lor 
\bigvee\limits_{\substack{j\\(j \neq r, j \neq i)}}
\neg X(r,i,j)
$
\item $
\bigl(\ \forall j \neq i: X(r,j,i)\ \bigr) \rightarrow A(r,i,n-1)
\ \ \ \ \ \Leftrightarrow \ \ \ \ \ 
\bigwedge\limits_{\substack{r,i\\(i \neq r)}}
A(r,i,n-1) \lor 
\bigvee\limits_{\substack{j\\(j \neq r, j \neq i)}}
\neg X(r,j,i)
$
\end{enumerate}
\normalsize

\item Connection between ``immediately after'' ($G$) and ``somewhere after'' ($X$):
$$
\forall r,i,j (i \neq r, j \neq r, i \neq j): G(r,i,j) \rightarrow X(r,i,j)
\ \ \ \ \ \Leftrightarrow \ \ \ \ \ 
\bigwedge\limits_{\substack{r,i,j\\(i \neq r, j \neq r, i \neq j)}}
\neg G(r,i,j) \lor X(r,i,j)
$$

\item \textit{(Consistency + Optimality)} 
Consistency was previously defined in the general properties: the relative ordering of any two lines $i$ and $j$ in row~$j$ must be consistent with how (pseudo)line intersections behave in the triangle formed by $i$, $j$, and $r$ (see Figure~\ref{fig:consistency}). Because in our case every finite non-overlapping segment must also be a side of a non-overlapping triangle, we enforce this ``optimality'' almost in the same way as we enforce consistency. The clauses are identical, except that for consistency we use $X$ (somewhere after), and for optimality we use $G$ (immediately after). The clauses below are for $G$, but must also be included identically for~$X$:

\resizebox{0.94\textwidth}{!}{
\begin{tabular}{ c c c }
\textit{Ordering} & \textit{First-order statement} & \textit{CNF clauses} \\

\hline

\raisebox{-.5\height}{\includegraphics[width=4cm]{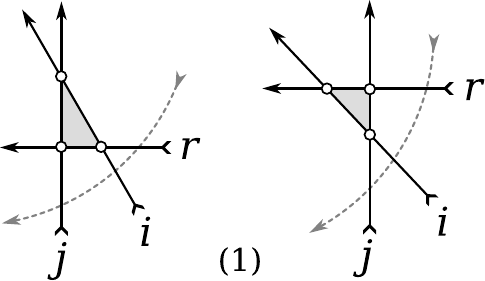}} & 
\begin{tabular}{c}
$(1)\ r < i < j\ :$ \\
$G(r,i,j) \leftrightarrow G(i,r,j) \land G(j,r,i))$ \\
$G(r,j,i) \leftrightarrow G(i,j,r) \land G(j,i,r))$
\end{tabular} & 
\begin{tabular}{c}
$\neg G(r,i,j) \lor G(i,r,j)$ \\
$\neg G(r,i,j) \lor G(j,r,i)$ \\
$G(r,i,j) \lor \neg G(i,r,j) \lor \neg G(j,r,i)$ \\
$\neg G(r,j,i) \lor G(i,j,r)$ \\
$\neg G(r,j,i) \lor G(j,i,r)$ \\
$G(r,j,i) \lor \neg G(i,j,r) \lor \neg G(j,i,r)$
\end{tabular} \\

\hdashline[0.3pt/0.8pt]

\raisebox{-.5\height}{\includegraphics[width=4cm]{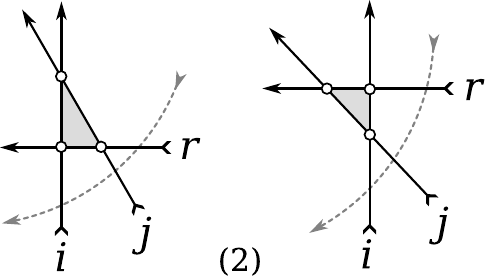}} & 
\begin{tabular}{c}
$(2)\ r < j < i\ :$ \\
$G(r,j,i) \leftrightarrow G(j,r,i) \land G(i,r,j))$ \\
$G(r,i,j) \leftrightarrow G(j,i,r) \land G(i,j,r))$
\end{tabular} & 
\begin{tabular}{c}
$\neg G(r,j,i) \lor G(j,r,i)$ \\
$\neg G(r,j,i) \lor G(i,r,j)$ \\
$G(r,j,i) \lor \neg G(j,r,i) \lor \neg G(i,r,j)$ \\
$\neg G(r,i,j) \lor G(j,i,r)$ \\
$\neg G(r,i,j) \lor G(i,j,r)$ \\
$G(r,i,j) \lor \neg G(j,i,r) \lor \neg G(i,j,r)$
\end{tabular} \\

\hdashline[0.3pt/0.8pt]

\raisebox{-.5\height}{\includegraphics[width=4cm]{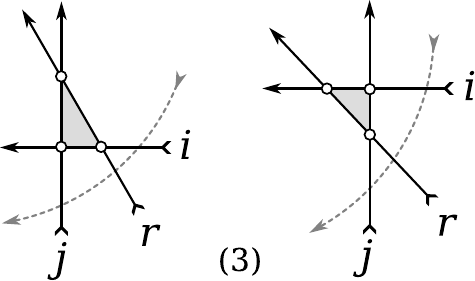}} & 
\begin{tabular}{c}
$(3)\ i < r < j\ :$ \\
$G(r,i,j) \leftrightarrow G(i,r,j) \land G(j,i,r))$ \\
$G(r,j,i) \leftrightarrow G(i,j,r) \land G(j,r,i))$
\end{tabular} & 
\begin{tabular}{c}
$\neg G(r,i,j) \lor G(i,r,j)$ \\
$\neg G(r,i,j) \lor G(j,i,r)$ \\
$G(r,i,j) \lor \neg G(i,r,j) \lor \neg G(j,i,r)$ \\
$\neg G(r,j,i) \lor G(i,j,r)$ \\
$\neg G(r,j,i) \lor G(j,r,i)$ \\
$G(r,j,i) \lor \neg G(i,j,r) \lor \neg G(j,r,i)$
\end{tabular} \\

\hdashline[0.3pt/0.8pt]

\raisebox{-.5\height}{\includegraphics[width=4cm]{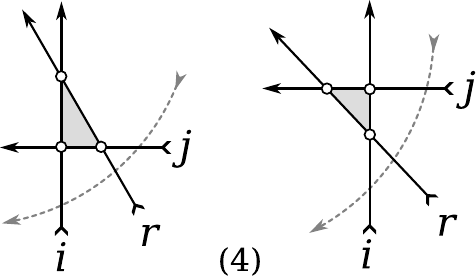}} & 
\begin{tabular}{c}
$(4)\ j < r < i\ :$ \\
$G(r,j,i) \leftrightarrow G(j,r,i) \land G(i,j,r))$ \\
$G(r,i,j) \leftrightarrow G(j,i,r) \land G(i,r,j))$
\end{tabular} & 
\begin{tabular}{c}
$\neg G(r,j,i) \lor G(j,r,i)$ \\
$\neg G(r,j,i) \lor G(i,j,r)$ \\
$G(r,j,i) \lor \neg G(j,r,i) \lor \neg G(i,j,r)$ \\
$\neg G(r,i,j) \lor G(j,i,r)$ \\
$\neg G(r,i,j) \lor G(i,r,j)$ \\
$G(r,i,j) \lor \neg G(j,i,r) \lor \neg G(i,r,j)$
\end{tabular} \\

\hdashline[0.3pt/0.8pt]

\raisebox{-.5\height}{\includegraphics[width=4cm]{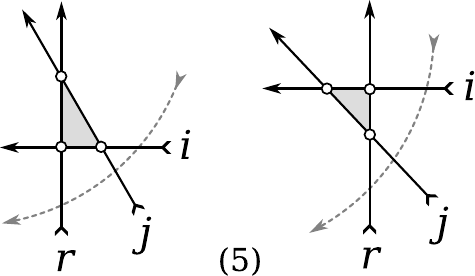}} & 
\begin{tabular}{c}
$(5)\ i < j < r\ :$ \\
$G(r,i,j) \leftrightarrow G(i,j,r) \land G(j,i,r))$ \\
$G(r,j,i) \leftrightarrow G(i,r,j) \land G(j,r,i))$
\end{tabular} & 
\begin{tabular}{c}
$\neg G(r,i,j) \lor G(i,j,r)$ \\
$\neg G(r,i,j) \lor G(j,i,r)$ \\
$G(r,i,j) \lor \neg G(i,j,r) \lor \neg G(j,i,r)$ \\
$\neg G(r,j,i) \lor G(i,r,j)$ \\
$\neg G(r,j,i) \lor G(j,r,i)$ \\
$G(r,j,i) \lor \neg G(i,r,j) \lor \neg G(j,r,i)$
\end{tabular} \\

\hdashline[0.3pt/0.8pt]

\raisebox{-.5\height}{\includegraphics[width=4cm]{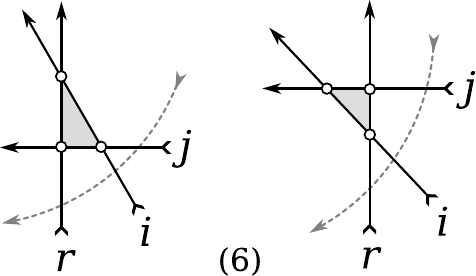}} & 
\begin{tabular}{c}
$(6)\ j < i < r\ :$ \\
$G(r,j,i) \leftrightarrow G(j,i,r) \land G(i,j,r))$ \\
$G(r,i,j) \leftrightarrow G(j,r,i) \land G(i,r,j))$
\end{tabular} & 
\begin{tabular}{c}
$\neg G(r,j,i) \lor G(j,i,r)$ \\
$\neg G(r,j,i) \lor G(i,j,r)$ \\
$G(r,j,i) \lor \neg G(j,i,r) \lor \neg G(i,j,r)$ \\
$\neg G(r,i,j) \lor G(j,r,i)$ \\
$\neg G(r,i,j) \lor G(i,r,j)$ \\
$G(r,i,j) \lor \neg G(j,r,i) \lor \neg G(i,r,j)$
\end{tabular} \\

\end{tabular} }

\item \textit{(Optional)} Mirror symmetry (mirrored across a line perpendicular to the first line):
\begin{enumerate}
\item (first line) $A(1,i,k) \rightarrow A(1, n-i+2, n-k)$
\item (other lines) $A(r,i,k) \rightarrow A(n-r+2, 1 + ((n-i+1) \bmod n), k)$
\end{enumerate}

\item \textit{(Optional)} Rotational symmetry: $A(r,i,k) \rightarrow A(sr,si,sk)$.

\item \textit{(Optional)} Let $\{t_{rc}\}$ be a table we wish to exclude. We add: $\bigvee\limits_{r,c} \neg A(r,t_{rc},c)$

\item \textit{(Optional)} 
By optimality, we expect each line to have every finite non-overlapping segment as a side of a non-overlapping triangle, but this is generally not the case for $n \bmod 6 \notin \{3,5\}$. We can generalize the model by allowing some lines to have ``missing triangles,'' meaning a line can have a finite non-overlapping segment that is not a side of a non-overlapping triangle. Let $L = \{l_1,l_2,\dots,l_m\}$ be a multiset of lines with missing triangles (a line included multiple times can have multiple missing triangles):
\begin{enumerate}
\item $\forall k: M_k(l_k,i,j) \rightarrow 
[\ \forall\ i',j'\ (i'\neq i \lor j'\neq j): \neg M_k(l_k,i',j')\ ]$
\item $\forall k: M_k(l_k,i,j) \rightarrow G(l_k,i,j)$
\item Add $M_k(l_k,i,j)$ to relevant optimality clauses: \\
$M_k(l_k,i,j) \lor [\ G(l_k,i,j) \leftrightarrow G(i,l_k,j) \land G(j,l_k,i)\ ],\ \ \dots$
\end{enumerate}

\end{enumerate}

\subsection{Statistics}

For a given $n$ and various options (such as symmetries or missing triangles), \texttt{kobon-cnf}~\cite{koboncnf} (specifically \texttt{gen.py}) generates a SAT problem in DIMACS CNF format and attempts to find a satisfying assignment using the Kissat~\cite{kissat} SAT solver. When successful, it produces a list of candidate tables. Another tool in the same repository (\texttt{fit.py}) helps to find the best straight-line arrangement corresponding to these candidate tables. Together, these tools can quickly reproduce many previously known Kobon arrangements and also generate new results (for example, maximal arrangements for $n=23$ and $n=27$, proving $N(23) = 161$ and $N(27) = 225$). See Table~\ref{tab:koboncnf-stats} for more detailed statistics, and Appendix~\ref{appendix:big_ar}.

\begin{table}
\begin{center}
\caption{Kobon-CNF Statistics\textsuperscript{\textdagger}}	
\label{tab:koboncnf-stats}
\
\begin{tabular}{r|ccc|cc|ll|ll|ll}
$n$ & \texttt{-M} & \texttt{-R} & \texttt{-L} 
& tabs & fit\textsuperscript{*} 
& $T_{\text{gen}}$\textsuperscript{**} & $T_{\text{fit}}$ 
& $\overline{T}_{\text{gen}}$ & $\overline{T}_{\text{fit}}$ 
& vars & clauses \\
\hline
3  & &   &     &  1 &  1 & 0.07 s  & 0.01 s  & 0.07 s  & 0.01 s  & 24    & 170\\
5  & &   &     &  1 &  1 & 0.08 s  & 0.015 s & 0.08 s  & 0.015 s & 200   & 2202\\
7  & &   & 3,6 &  1 &  1 & 0.3 s   & 0.14 s  & 0.3 s   & 0.14 s  & 732   & 11896\\
9  & &   &     &  1 &  1 & 0.54 s  & 0.35 s  & 0.54 s  & 0.35 s  & 1584  & 30102\\
11 & &   &     &  0 & -- & 1.67 s  & --      & --      & --      & 3080  & 70840\\
13 & &   & 6,9 &  3 &  3 & 84.7 s  & 4.83 s  & 28.23 s & 1.61 s  & 5568  & 178162\\
\hdashline[0.3pt/0.8pt]
15 & &   &     &  4 &  4 & 140.0 s & 18.57 s & 35.0 s  & 4.64 s  & 8400  & 260592\\
15 & & 3 &     &  0 & -- & 8.42 s  & --      & --      & --      & 8400  & 266280\\
15 & & 5 &     &  1 &  1 & 6.38 s  & 0.24 s  & 6.38 s  & 0.24    & 8400  & 272172\\
\hdashline[0.3pt/0.8pt]
17 & &   &     & 10 & 10 & 1307 s  & 195.5 s & 130.7 s & 19.5 s  & 12512 & 438430\\
19 & \checkmark & & 3,18 & 8 & 8 & 733.4 s & 262.2 s & 91.7 s & 32.8 s & 18396 & 887310 \\
\hdashline[0.3pt/0.8pt]
21 & & 3 &     & 10 &  4 & 620.1 s & 269.3 s & 62.0 s  & 26.9 s & 24360 & 1064560\\
21 & & 7 &     &  0 & -- & 35.8 s  & --      & --      & --     & 24360 & 1097880\\
21 & \checkmark & & & 4 & 4 & 1216 s & 237.9 s & 304 s & 59.4 s & 24360 & 1056048\\
21 & &   &     & 236 & 85 & 44h 16m & 10 h & 675 s  & 152 s & 24360 & 1066576\\
\hdashline[0.3pt/0.8pt]
23 & & & & $\geq65$ & $\geq6$ & 9h 25m & 2h 15h & 521 s & 125 s & 32384 & $\geq1522048$\\
23 & \checkmark & & & 0 & -- & 213 s & -- & -- & -- & 32384 & 1533180\\
\hdashline[0.3pt/0.8pt]
27 & & 3 & & 86 & 5 & 8 h & 35 m & 337 s & 24.4 s & 53352 & 2973960\\
27 & \checkmark & & & 0 & -- & 307 s & -- & -- & -- & 53352 & 2952612\\
27 & & 9 & & 0 & -- & 258 s & -- & -- & -- & 53352 & 3080376\\
\hdashline[0.3pt/0.8pt]
35 & \checkmark & & & 0 & -- & 1h 20m & -- & -- & -- & 119000 & 8489612\\
35 & & 7 & & 0 & -- & 1h 36m & -- & -- & -- & 119000 & 8691760\\
\hdashline[0.3pt/0.8pt]
39 & & 13 & & 0 & -- & 1h 28m & -- & -- & -- & 165984 & 13788528
\end{tabular}
\end{center}
\begin{tabular}{ll}
$n$ & -- number of lines, \\
\texttt{-M} & -- mirror symmetry, \\
\texttt{-R} & -- rotational symmetry, \\
\texttt{-L} & -- lines with a missing triangle, \\
$T_{\text{gen}}$ & -- time to generate the tables; 
$\overline{T}_{\text{gen}} = T_{\text{gen}} / tabs$, \\
$T_{\text{fit}}$ & -- time to fit straight lines;
$\overline{T}_{\text{fit}} = T_{\text{fit}} / tabs$, \\
$\geq k$ & -- search not exhaustive; more arrangements may be found with \texttt{gen.py}.
\end{tabular}

\smallskip
\noindent\rule{3cm}{0.4pt}
\smallskip

\begin{minipage}{\linewidth}
\footnotesize
{\textsuperscript{\textdagger}
Experiments performed on ASUS GL704GW, Kissat v4.0.2, SciPy v1.7.3, Python 3.7.4.}

* Column ``fit'' shows the number of unique and successfully straightened candidates using default parameters (except for $n=27$ where we use \texttt{main\_coefficient=0.5}).

** Including the time to prove that no other tables exist.
\end{minipage}
\end{table}

\begin{figure}
\centering
\includegraphics[width=0.85\textwidth]{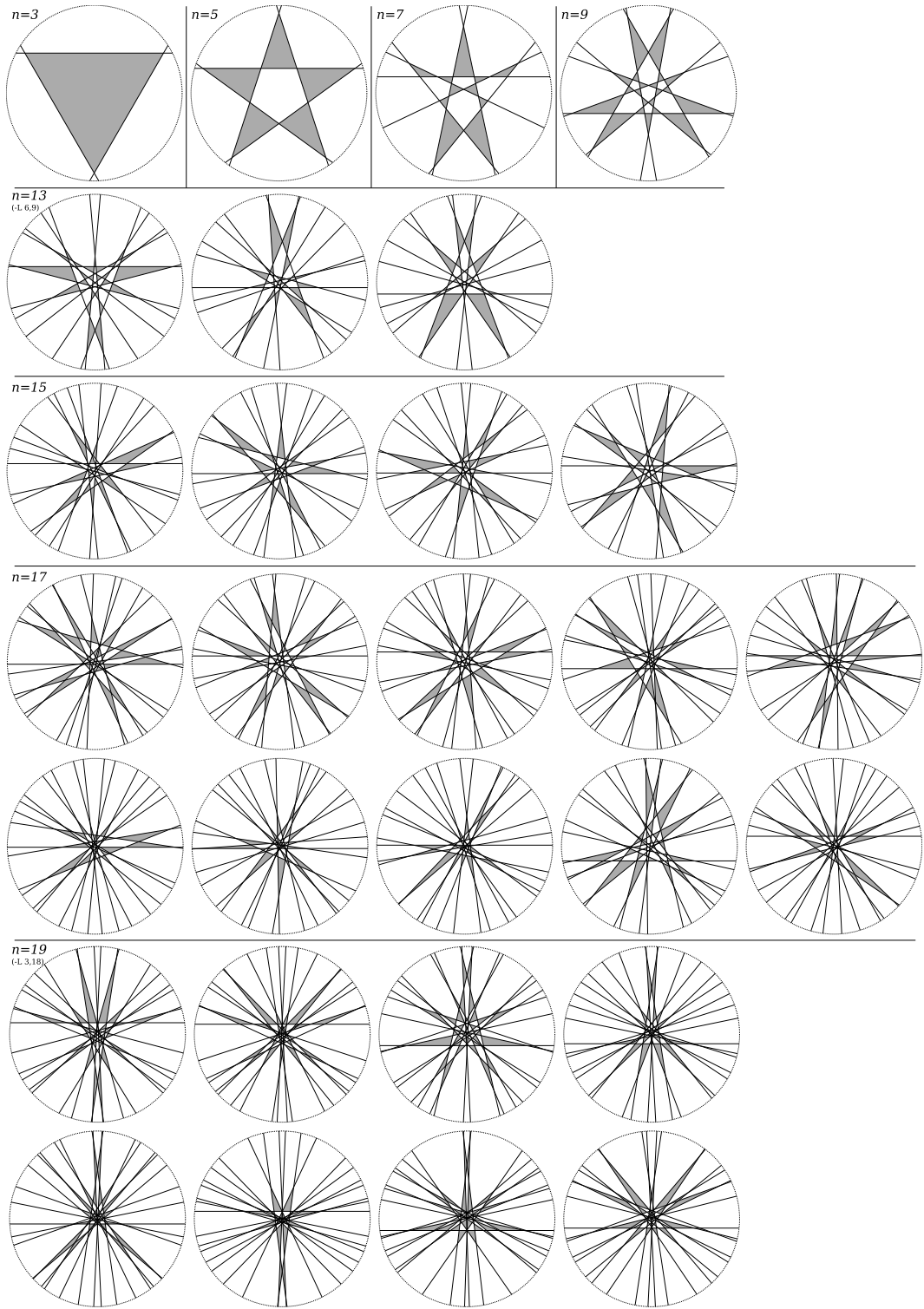} \\
\caption{Optimal arrangements found by \texttt{kobon-cnf} for $n\in\{3,5,7,9,13,15,17,19\}$. For $n\in\{13,19\}$, other arrangements could exist (for example, with different lines having missing triangles).}
\label{fig:kobon-cnf-all-01}
\end{figure}

\section{Heuristic Straightening}

In this section, we give a brief overview of the heuristic straightening method used to find balanced arrangements of straight lines that correspond to a given table, including cases with symmetries and other properties. It is a surprisingly simple approach --  essentially a constrained minimization of a cleverly constructed function, using \texttt{SciPy} and \texttt{NumPy}. Despite its simplicity, it works well for many complex cases. Every single arrangement image shown in this paper was created using this approach. See the documented code of \texttt{solver(..)} in~\cite{lineorder} for a working realization.

\subsection{Unknown Variables}

Each unknown line is represented by two variables: $a_i$ (its angle) and $C_i$ (its signed distance to the origin). Therefore, the equation of the $i$-th line is
$$x \cos a_i + y \sin a_i + C_i = 0, \quad i = 1, 2, \dots, n.$$

For a given line $i$, two consecutive entries $j$ and $k$ in the $i$-th table row specify that line $j$ must intersect line $i$ before line $k$ intersects line $i$. With $a_i, a_j, a_k$ and $C_i, C_j, C_k$, this condition can be expressed using an inequality:
\begin{center}
\begin{tabular}{l}
$S(i,j,k) \cdot F(i,j,k) < 0$, where \\ 
$F(i,j,k) = 
    C_i \cdot \sin(a_k - a_j) + 
    C_j \cdot \sin(a_i - a_k) + 
    C_k \cdot \sin(a_j - a_i)$
\end{tabular}
\end{center}
At a multi-line intersection point, where three or more lines intersect at the same point, these inequalities become equalities among the lines involved.

\begin{samepage}
$S(i,j,k)$ is a special ``sign-correction'' function, defined in Python code:
\begin{center}
\small
\begin{BVerbatim}
def cmp_func(row, l1, l2):
    if l1 < row: l1 += n
    if l2 < row: l2 += n
    return l1 < l2

def S(row, l1, l2):
    s = 1.0 if l1 > l2 else -1.0
    s *= -1.0 if cmp_func(row, l1, l2) else 1.0
    return s
\end{BVerbatim}
\normalsize
\end{center}
\end{samepage}

\begin{samepage}
The other limitations are for line angles $a_i$ and signed distances $C_i$:
\begin{center}
\begin{tabular}{c}
$a_i > a_{i+1}$, \ \\ 
$-\pi < a_i < 0$, \ \\ 
$A_{min} < |a_i - a_{i+1}| < A_{max}$, \ \\ 
$C_{min} < C_i < C_{max}$.
\end{tabular}
\end{center}
\end{samepage}

In the summary, for $n$ lines, we have no more than $2n$ unknowns. The number of unknowns can be reduced if it is known that the arrangement has symmetries. E.g. with \texttt{rot3}-symmetry, only $1/3$ of variables are independent.

\subsection{Target Function}

At its core, the solver works by minimizing a specially designed target function. Smaller values of this function correspond to smaller violations of the imposed conditions. The conditions for the inequalities are codified using:
$$
\text{LessThan}(x,v) = 
\begin{cases}
 0, & \text{ if } \ x < v \\
 (x-v)^2, & \text{ if } \ x >= v
\end{cases}
$$
The target function is defined as a sum of \texttt{LessThan} functions with coefficients:
\begin{enumerate}
\item Regular cross-points (neighboring $l_1$ and $l_2$ in row $r$):
$$-\text{INEQ}_{\max} < F(r, l_1, l_2) < -\text{INEQ}_{\epsilon} < 0$$
\item Multiline cross-points (a group $l_1,l_2,\dots,l_k$ crossing $r$ at the same intersection):
$$|F(r, l_i, l_{i+1})| \leq 0$$
\item Angle constraints: $a_i < a_0$
\item The angle between two consecutive lines is in between $A_{\min}$ and $A_{\max}$:
$$A_{\min} < a_{i} - a_{i+1} < A_{\max}$$
\item The angle between $a_1$ and $a_n$ is greater than $A_{\min}$:
$$a_1 - \pi < a_n - A_{\min}$$
\end{enumerate}
Some conditions are encoded directly as bounds on the variables:
\begin{enumerate}
\item $0 \leq -a_i \leq \pi/S$,\ \ where $S$ is the rotational symmetry (1, 3, \dots).
\item $-100 \leq C_i \leq 100$
\end{enumerate}

\subsection{Initial Guess}

The initial guess $x_0$ for the \texttt{minimize} function is defined as:
$$a_i = -\pi/(2n) - \pi(i-1)/n$$
$$C_i = 0.1\cdot(-1)^{i\bmod 2}$$

\subsection{Minimization}

With \texttt{scipy.optimize.minimize}, we attempt to find the minimum of the target function. Any set of variables (within the boundaries) that produce a value of zero is considered an ideal solution. However, in practice, even nonzero values often correspond to perfectly acceptable solutions for the overall problem of straightening an arrangement.

\subsection{Forcing Additional Properties}

The method is flexible enough to enforce additional properties and symmetries in the line arrangements. For example, it is possible to force the lines to take the form $y = m_i(x-a_i)$ where $a_i$ are either very small values or $\tan\bigl(\pm i \cdot \pi/(n-1)\bigr)$, which is a form suitable for Proposition~3.1 in~\cite{bartholdi2007}. See Appendix~\ref{appendix:tan_form} for the resulting straight-line arrangements.

\section{Conclusion}
The approach described in this paper has proven productive, enabling the discovery of new maximal Kobon arrangements (for $n\in\{23,27\}$). For $n=11$ we confirm that a simple perfect arrangement with 33 triangles cannot exist, even among pseudolines. For $n\in\{3,5,9,15,17\}$ we enumerate all possible simple perfect Kobon arrangements.

One possible direction for further research is to generate maximal tables for even $n$ and attempt to straighten them using our heuristic straightening method via \texttt{LineOrder}~\cite{lineorder}. These tables may also include parallel lines and multi-line cross points, which are underutilized in this research. 
Another interesting topic is the analysis of when and why heuristic straightening works well and when it does not. This is connected to the question of how to determine in advance whether an arrangement can be straightened at all (which is generally a hard problem~\cite{Fukuda2013}), and how to modify the CNF model to exclude many unstraightenable arrangements.
Finally, it is still possible to find new optimal arrangements using \texttt{Kobon-CNF}~\cite{koboncnf} through a trial-and-error strategy.

\section{Acknowledgments}

\texttt{LineOrder}\cite{lineorder} was created during the collaboration with {Kyle Wood} on Kobon Triangles. Special thanks to Kyle for providing complex test cases featuring mirror symmetry. \texttt{LineOrder} is written in Python~3.8 and uses \texttt{SciPy}\cite{2020SciPy-NMeth}, \texttt{NumPy}\cite{harris2020array} and \texttt{Matplotlib}\cite{Hunter:2007} libraries. \texttt{Kobon-CNF}\cite{koboncnf} is written in Python~3.8 and uses \texttt{Kissat}~SAT~solver\cite{kissat}. 

\appendix

\section{Special Arrangements}
\label{appendix:tan_form}
This appendix presents arrangements in the form $y = m_i(x-a_i),$ where the $a_i$ are either very small values or of the form $\tan(\pm i\cdot\pi/(n-1))$, making them suitable for Proposition~3.1 in \cite{bartholdi2007}. These arrangements were generated using \texttt{LineOrder}\cite{lineorder}, a tool that applies the heuristic straightening approach described in this paper:

\subsection{Kobon 5-Lines (5 Triangles, Mirror Symmetry)}
\begin{samepage}
\begin{center}
\begin{tabular}{rc}
\begin{tabular}{ll}
$\epsilon = 1 / 10 < 1 / 8$ & \\
\multicolumn{2}{l}{$y = m_i(x - a_i)$, $i = 0..4:$} \tabularnewline
\hdashline[0.3pt/0.8pt]
$m_0 = 0$,           & $a_0 = 0$\\
$m_1 = -1.3763819$,  & $a_1 = \tan(-1\pi/4)$\\
$m_2 = -19.9833325$, & $a_2 = \epsilon$\\
$m_3 = 19.9833325$,  & $a_3 = -\epsilon$\\
$m_4 = 1.3763819$,   & $a_4 = \tan(1\pi/4)$
\end{tabular} 
& \raisebox{-.5\height}{\includegraphics[width=3cm]{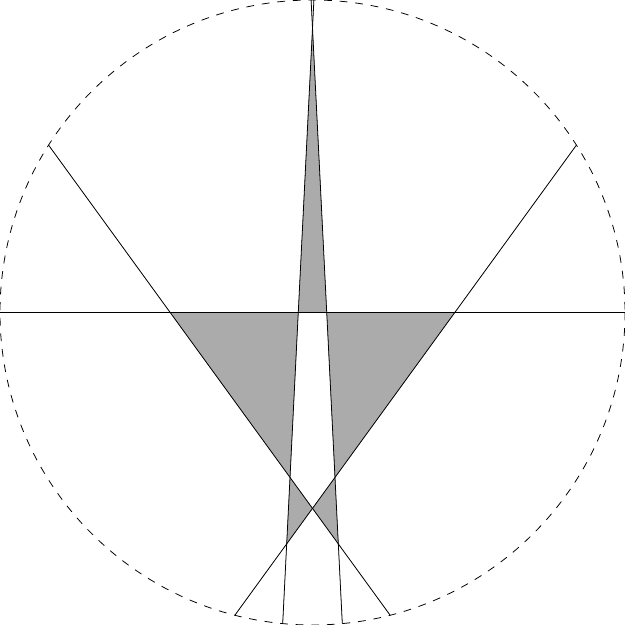}}
\end{tabular}
\end{center}
\normalsize
\end{samepage}

\subsection{Kobon 7-Lines (11 Triangles, Mirror Symmetry)}
\begin{samepage}
\small
\begin{center}
\begin{tabular}{rc}
\begin{tabular}{ll}
$\epsilon = 1 / 14 < 1 / 12$ & \\
\multicolumn{2}{l}{$y = m_i(x - a_i)$, $i = 0..6:$} \tabularnewline
\hdashline[0.3pt/0.8pt]
$m_0 = 0$,           & $a_0 = 0$\\
$m_1 = -0.7974734$,  & $a_1 = \tan(-1\pi/6)$\\
$m_2 = -2.0765214$,  & $a_2 = \tan(-2\pi/6)$\\
$m_3 = -19.9833325$, & $a_3 = \epsilon$\\
$m_4 = 19.9833325$,  & $a_4 = -\epsilon$\\
$m_5 = 2.0765214$,   & $a_5 = \tan(2\pi/6)$\\
$m_6 = 0.7974734$,   & $a_6 = \tan(1\pi/6)$
\end{tabular} 
& \raisebox{-.5\height}{\includegraphics[width=4cm]{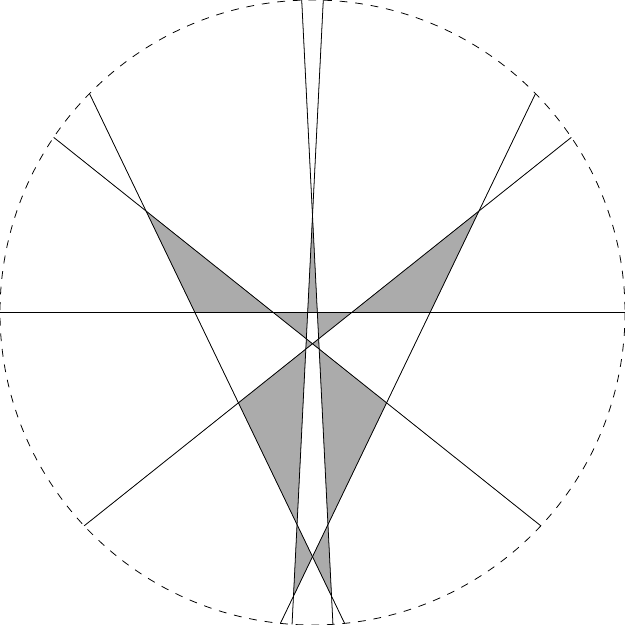}}
\end{tabular}
\end{center}
\normalsize
\end{samepage}

\subsection{Kobon 9-Lines (21 Triangles, Mirror Symmetry)}
\begin{samepage}
\small
\begin{center}
\begin{tabular}{rc}
\begin{tabular}{ll}
$\epsilon = 1 / 18 < 1 / 16$ & \\
\multicolumn{2}{l}{$y = m_i(x - a_i)$, $i = 0..8:$} \tabularnewline
\hdashline[0.3pt/0.8pt]
$m_0 = 0$,          & $a_0 = 0$\\
$m_1 = -0.5773503$, & $a_1 = \tan(-3\pi/8)$\\
$m_2 = -1.1917536$, & $a_2 = \tan(-1\pi/8)$\\
$m_3 = -2.7474774$, & $a_3 = \tan(-2\pi/8)$\\
$m_4 = -19.9833325$,& $a_4 = \epsilon$\\
$m_5 = 19.9833325$, & $a_5 = -\epsilon$\\
$m_6 = 2.7474774$,  & $a_6 = \tan(2\pi/8)$\\
$m_7 = 1.1917536$,  & $a_7 = \tan(1\pi/8)$\\
$m_8 = 0.5773503$,  & $a_8 = \tan(3\pi/8)$
\end{tabular} 
& \raisebox{-.5\height}{\includegraphics[width=4.5cm]{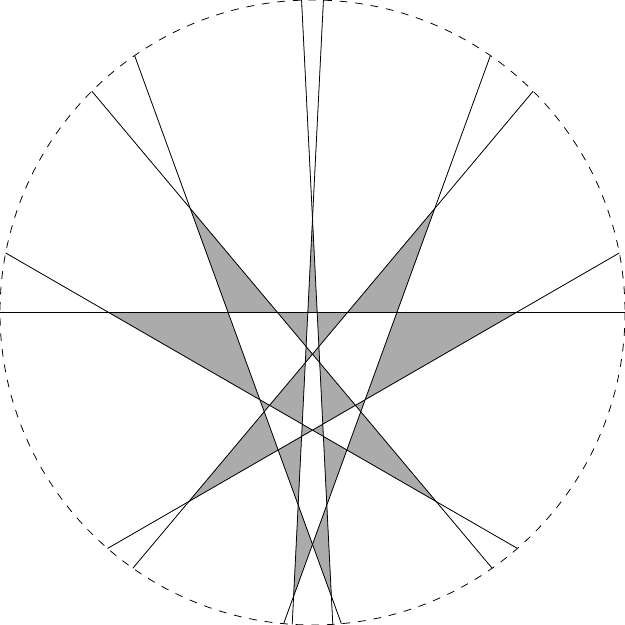}}
\end{tabular}
\end{center}
\normalsize
\end{samepage}

\subsection{11-Lines (32 Triangles)}
\begin{samepage}
Based on solution by Honma, not optimal (32/33 triangles):
\small
\begin{center}
\begin{tabular}{rc}
\begin{tabular}{ll}
$\epsilon = 1 / 22 < 1 / 20$ & \\
\multicolumn{2}{l}{$y = m_i(x - a_i)$, $i = 0..10:$} \tabularnewline
\hdashline[0.3pt/0.8pt]
$m_0 = 0$,           & $a_0 = 0$\\
$m_1 = -0.2615465$,  & $a_1 = \tan(-4\pi/10)$\\
$m_2 = -1.0298714$,  & $a_2 = \tan(1\pi/10)$\\
$m_3 = -2.0130975$,  & $a_3 = \tan(-2\pi/10)$\\
$m_4 = -4.3446427$,  & $a_4 = -\epsilon$\\
$m_5 = -25.5645486$, & $a_5 = \tan(-3\pi/10)$\\
$m_6 = 2.7017707$,   & $a_6 = \tan(-1\pi/10)$\\
$m_7 = 1.4888582$,   & $a_7 = \epsilon$\\
$m_8 = 1.0746682$,   & $a_8 = \tan(2\pi/10)$\\
$m_9 = 0.718313$,    & $a_9 = \tan(3\pi/10)$\\
$m_{10} = 0.1670701$,  & $a_{10} = \tan(4\pi/10)$
\end{tabular} 
& \raisebox{-.5\height}{\includegraphics[width=5.5cm]{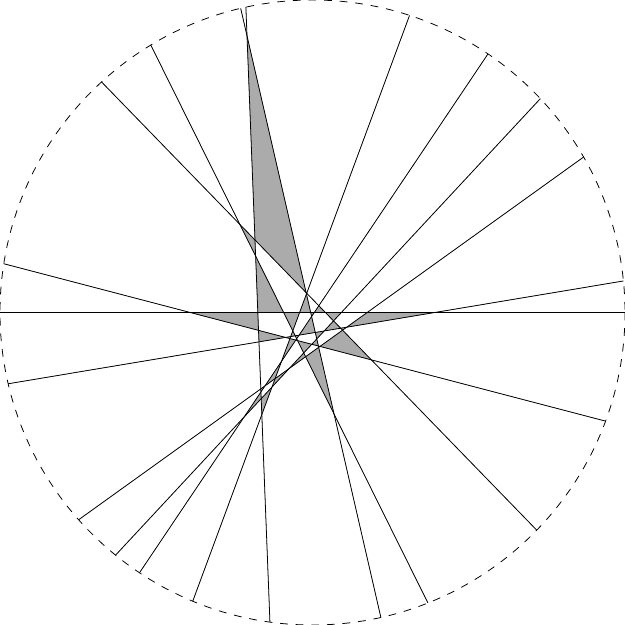}}
\end{tabular}
\end{center}
\normalsize
\end{samepage}

\subsection{Kobon 13-Lines (47 Triangles, Mirror Symmetry)}
\begin{samepage}
Based on solution by Kabanovitch:
\small
\begin{center}
\begin{tabular}{rc}
\begin{tabular}{ll}
$\epsilon = 1 / 26 < 1 / 24$ & \\
\multicolumn{2}{l}{$y = m_i(x - a_i)$, $i = 0..12:$} \tabularnewline
\hdashline[0.3pt/0.8pt]
$m_0 = 0$,           & $a_0 = 0$\\
$m_1 = -0.2885669$, & $a_1 = \tan(-5\pi/12)$\\
$m_2 = -0.7122015$, & $a_2 = \tan(-1\pi/12)$\\
$m_3 = -0.8749263$, & $a_3 = \tan(-3\pi/12)$\\
$m_4 = -3.6695876$, & $a_4 = \tan(-2\pi/12)$\\
$m_5 = -6.0310043$, & $a_5 = \tan(-4\pi/12)$\\
$m_6 = -15.8845208$,& $a_6 = \epsilon$\\
$m_7 = 15.8845208$, & $a_7 = -\epsilon$\\
$m_8 = 6.0310043$,  & $a_8 = \tan(4\pi/12)$\\
$m_9 = 3.6695876$,  & $a_9 = \tan(2\pi/12)$\\
$m_{10} = 0.8749263$, & $a_{10} = \tan(3\pi/12)$\\
$m_{11} = 0.7122015$, & $a_{11} = \tan(1\pi/12)$\\
$m_{12} = 0.2885669$, & $a_{12} = \tan(5\pi/12)$
\end{tabular} 
& \raisebox{-.5\height}{\includegraphics[width=6.5cm]{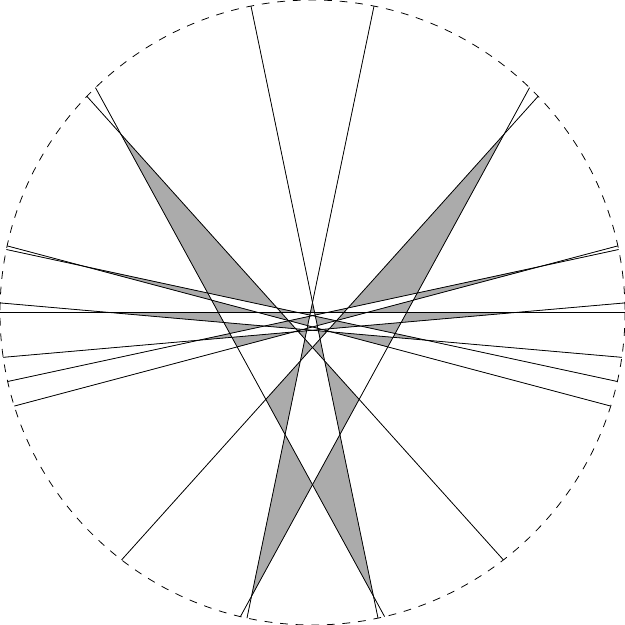}}
\end{tabular}
\end{center}
\normalsize
\end{samepage}

\subsection{Kobon 17-Lines (85 Triangles)}
\begin{samepage}
Based on solution by Johannes Bader:
\small
\begin{center}
\begin{tabular}{rc}
\begin{tabular}{ll}
$\epsilon = 1 / 34 < 1 / 32$ & \\
\multicolumn{2}{l}{$y = m_i(x - a_i)$, $i = 0..16:$} \tabularnewline
\hdashline[0.3pt/0.8pt]
$m_0 = 0$,           & $a_0 = 0$\\
$m_1 = -0.3166485$, & $a_1 = \tan(7\pi/16)$\\
$m_2 = -0.4791101$, & $a_2 = \tan(-4\pi/16)$\\
$m_3 = -0.8204082$, & $a_3 = \tan(3\pi/16)$\\
$m_4 = -1.0031728$, & $a_4 = \tan(-6\pi/16)$\\
$m_5 = -1.6804231$, & $a_5 = \tan(1\pi/16)$\\
$m_6 = -2.3663988$, & $a_6 = \tan(-2\pi/16)$\\
$m_7 = -5.3495275$, & $a_7 = \tan(5\pi/16)$\\
$m_8 = 26.7274944$, & $a_8 = -\epsilon$\\
$m_9 = 5.4191716$,  & $a_9 = \tan(2\pi/16)$\\
$m_{10} = 2.378765$,  & $a_{10} = \tan(-5\pi/16)$\\
$m_{11} = 1.5246315$, & $a_{11} = \epsilon$\\
$m_{12} = 1.1361431$, & $a_{12} = \tan(-3\pi/16)$\\
$m_{13} = 0.9298155$, & $a_{13} = \tan(4\pi/16)$\\
$m_{14} = 0.4979408$, & $a_{14} = \tan(-1\pi/16)$\\
$m_{15} = 0.241266$,  & $a_{15} = \tan(6\pi/16)$\\
$m_{16} = 0.1047156$, & $a_{16} = \tan(-7\pi/16)$
\end{tabular} 
& \raisebox{-.5\height}{\includegraphics[width=9.0cm]{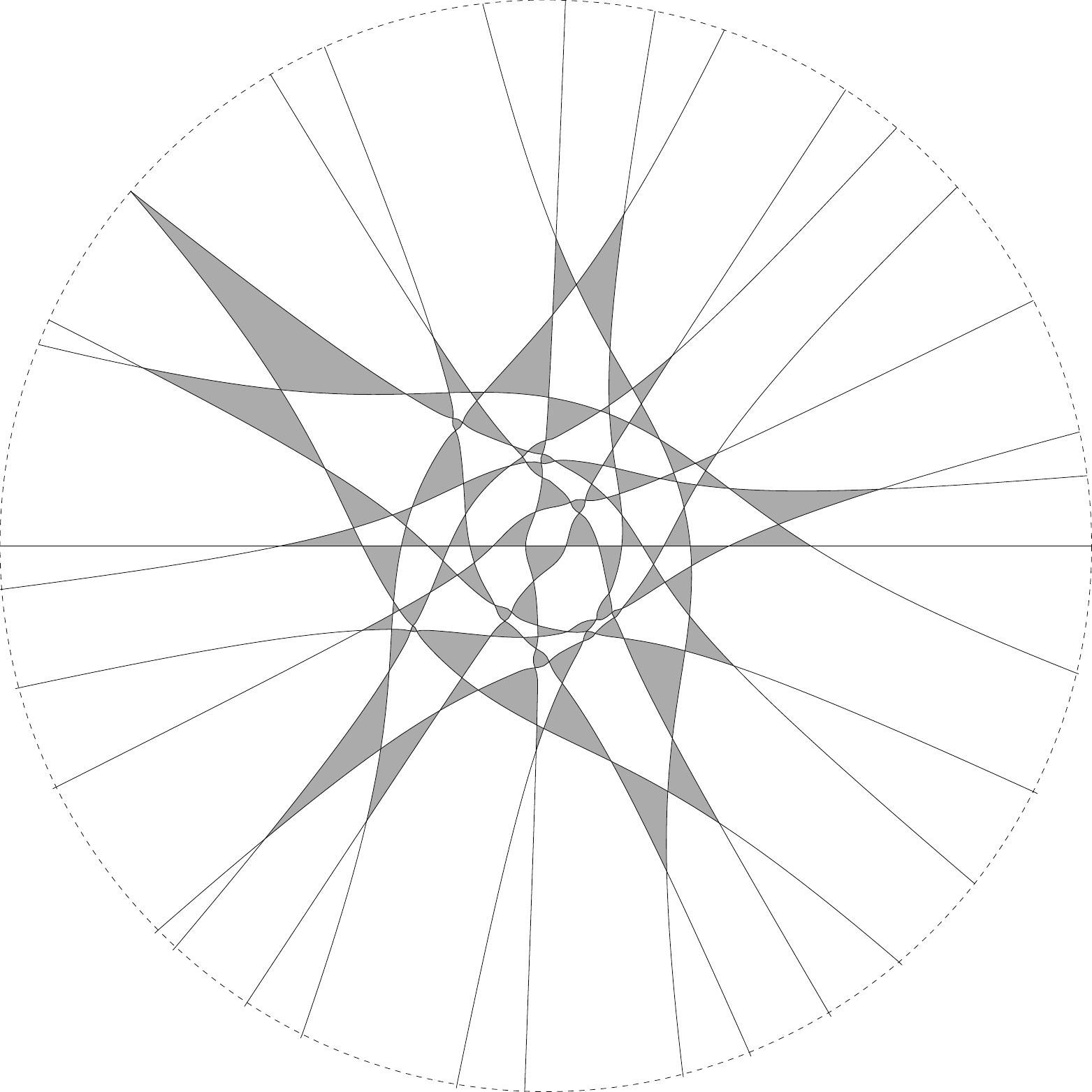}}
\end{tabular}
\end{center}
\normalsize
\end{samepage}
\pagebreak

\section{On Presentation}
\label{appendix:presentation}
Producing images that clearly communicate results is an important part of scientific research. As $n$ grows larger, Kobon arrangements become increasingly difficult to present visually. Drawings with straight lines are useful to demonstrate that an arrangement is stretchable, but for large $n$, such drawings become problematic. 

The first issue is that large arrangements tend to contain very small triangles near the center and much larger, elongated triangles toward the outside. The second issue is that the size of neighboring triangles can vary drastically, often by an order of magnitude. These issues make the resulting images less suitable as visual illustrations\cite{rote2015presentation}, especially for print.

\begin{figure}
\begin{center}
\includegraphics[width=11.5cm]{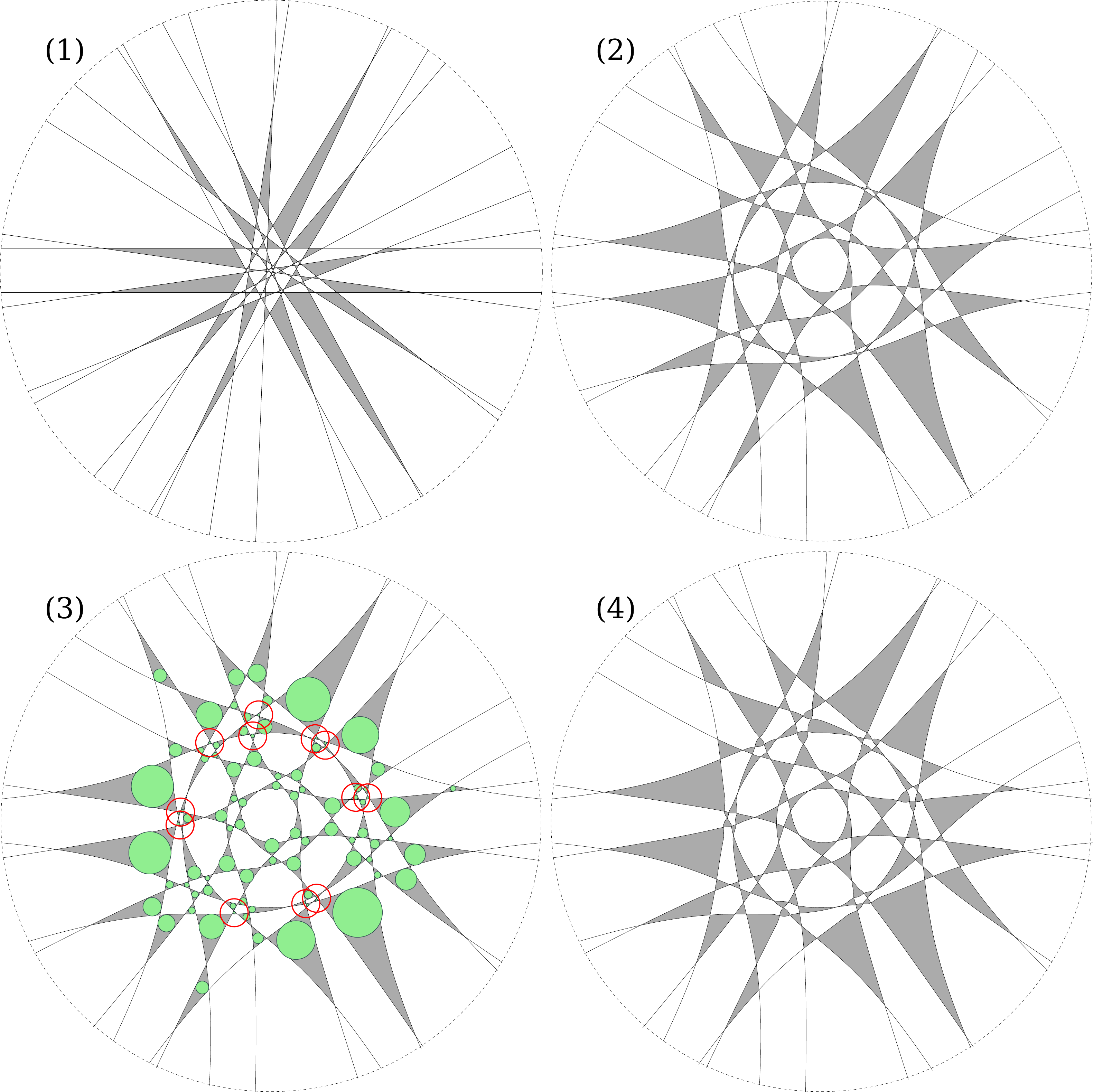}
\end{center}
\caption{These images show how the described techniques can be used to improve visual presentation. (1)~Straight-line drawing. (2)~Fisheye projection applied. (3)~Detection of small problematic triangles. (4)~Local zoom-in effect applied. (18-line arrangement with 93 triangles by Johannes Bader\cite{baderj_other}).}
\label{fig:better_drawing}
\end{figure}

The first issue can be addressed by applying a fisheye projection to reduce the size difference between medial and lateral triangles. This also shortens long lateral triangles (see Figure~\ref{fig:better_drawing}(1)-(2)). The second issue can be partially resolved by locally applying fisheye projection to clusters of very small triangles (see Figure~\ref{fig:better_drawing}(3)-(4)).	

For more technical details please visit \cite[Gallery~\#3]{lineorder} and \texttt{lineorder/draw\_lines.py}.

\pagebreak
\section{Big Arrangements}
\label{appendix:big_ar}
This section includes several new large arrangements found using the approaches described in this paper. The corresponding tables are provided in the last subsection of this appendix. More information is available on the ``Gallery'' page of~\cite{lineorder}, including images without fisheye projection.

\subsection{23-Line Solution (161 Triangles)}
A 23-line solution with 161 triangles, achieving the upper bound and proving that $N(23) = 161$:

\begin{center}
\includegraphics[width=\textwidth]{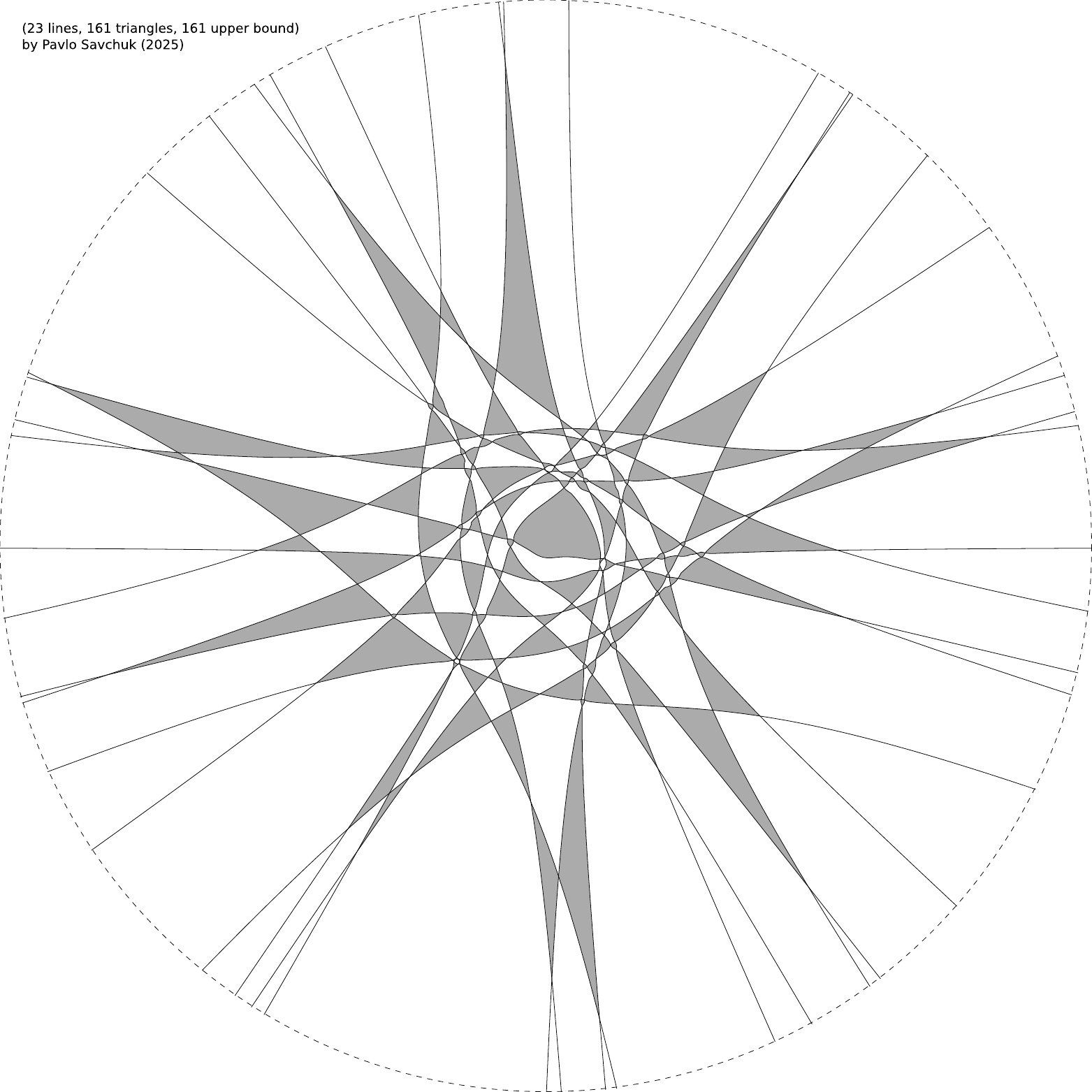}
\end{center}

\subsection{24-Line Solution (172 Triangles)}
Following the examples of the 16-line solution by Bader\cite{baderj_other} and the 20-line solution by Wood\cite[A006066]{oeis}, this arrangement was constructed by adding an additional line between the first and second lines of the 23-line arrangement from the previous subsection. In a similar way, it is possible to construct 22-line and 28-line arrangements based on the 21- and 27-line solutions, but as in this case, they also lack a triangle to meet the current best upper bounds. Such even tables can be auto-generated from odd tables using:
\begin{verbatim}
lineorder.add_1_2_line(table, ...)
\end{verbatim}

\begin{center}
\includegraphics[width=\textwidth]{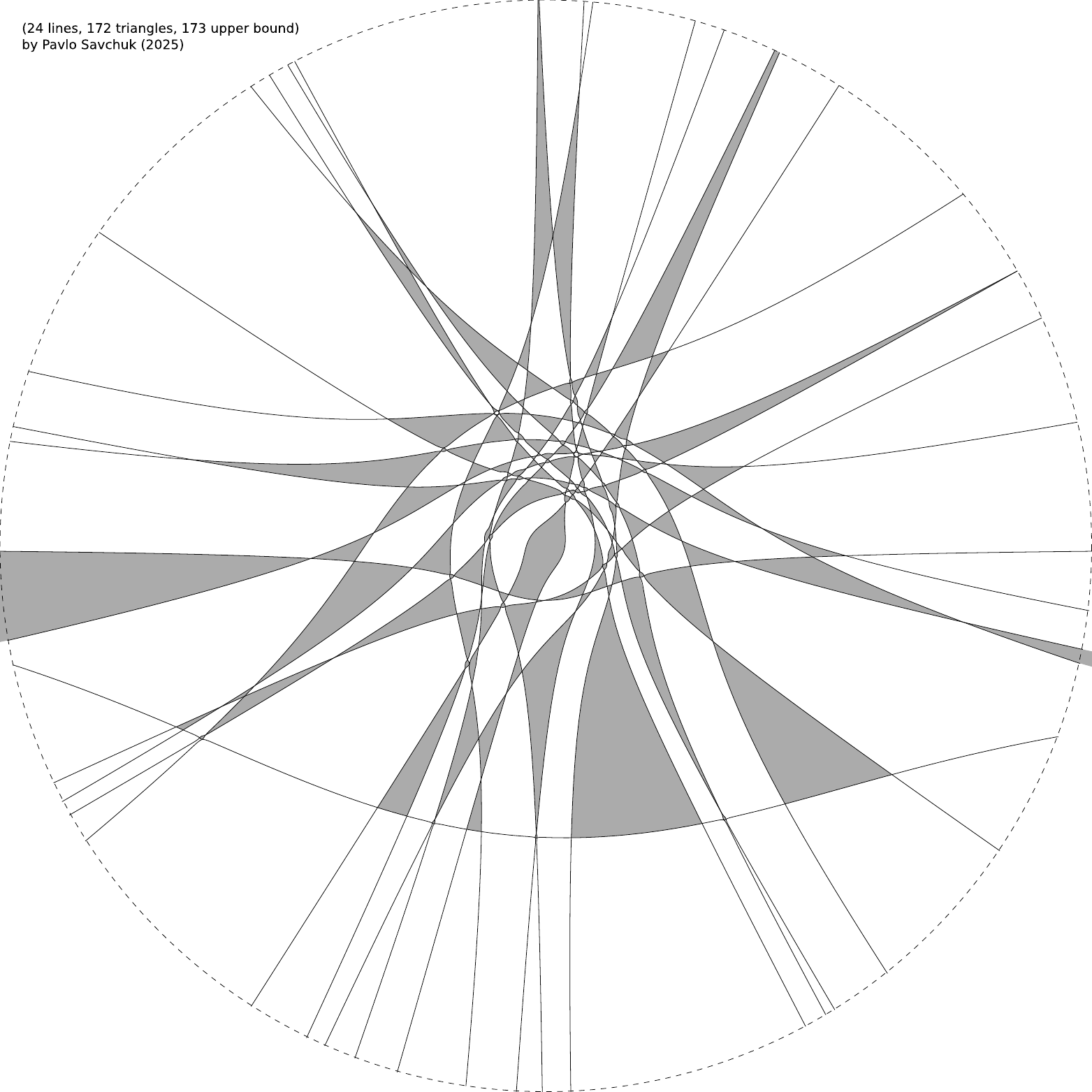}
\end{center}

\subsection{27-Line Solution \#1 (225 Triangles)}
\begin{samepage}
A 27-line solution \#1 with 225 triangles, achieving the upper bound and proving that $N(27) = 225$:

\vspace{0.5cm}
\begin{center}
\includegraphics[width=\textwidth]{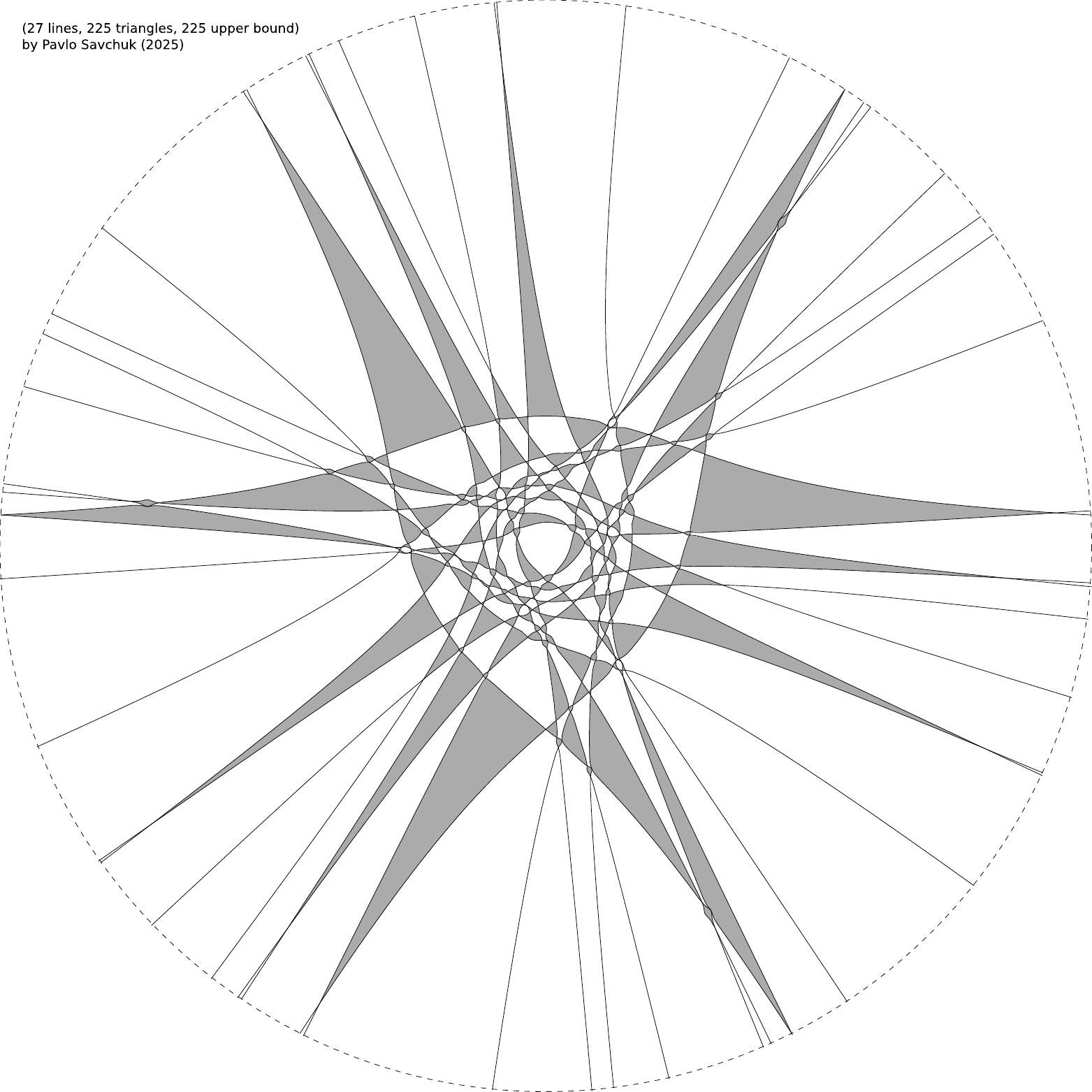}
\end{center}
\end{samepage}

\subsection{27-Line Solution \#2 (225 Triangles)}
\begin{samepage}
Another, more compact 27-line solution with 225 triangles, achieving the upper bound and proving that $N(27)=225$:

\begin{center}
\includegraphics[width=\textwidth]{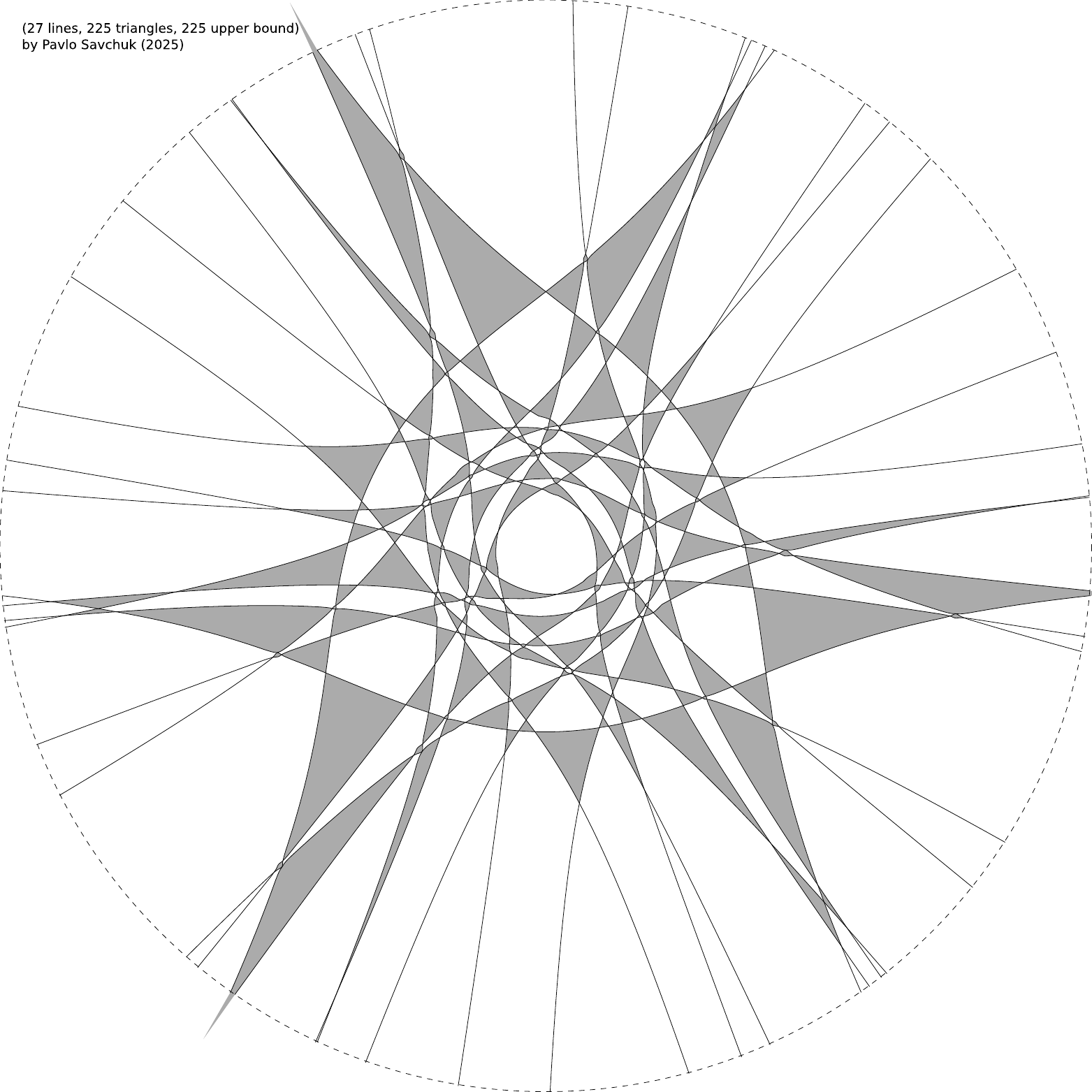}
\end{center}
\end{samepage}

\subsection{Corresponding Tables}
\begin{samepage}
These and many more examples can be found on the ``Gallery'' page of~\cite{lineorder}.\begin{center}
\begin{adjustbox}{max width=\linewidth}
\begin{tabular}{cc}
23-Line Solution (161 Triangles): &
24-Line Solution (172 Triangles): \\
\tiny
\begin{BVerbatim}
[[2,20,4,18,8,22,6,14,3,12,10,16,7,17,9,15,13,19,11,21,5,23],
 [1,20,22,18,21,16,19,14,17,8,15,12,23,10,13,6,9,7,11,4,5,3],
 [4,20,8,18,6,22,14,1,12,16,10,17,7,15,9,19,13,21,11,23,5,2],
 [3,20,1,18,22,8,16,14,21,12,17,6,15,10,19,7,13,9,23,11,2,5],
 [6,8,7,10,9,14,12,18,13,16,11,17,15,20,19,22,21,1,23,3,2,4],
 [5,8,20,18,3,22,1,14,16,12,21,17,4,15,19,10,23,13,2,9,11,7],
 [8,5,10,18,14,20,12,22,16,1,17,3,15,21,19,4,13,23,9,2,11,6],
 [7,5,6,20,3,18,1,22,4,16,21,14,19,17,2,15,23,12,13,10,11,9],
 [10,5,14,18,12,20,16,22,17,1,15,3,19,21,13,4,23,7,2,6,11,8],
 [9,5,7,18,20,14,22,12,1,16,3,17,21,15,4,19,6,23,2,13,8,11],
 [12,14,13,18,16,5,17,20,15,22,19,1,21,3,23,4,2,7,6,9,8,10],
 [11,14,5,18,9,20,7,22,10,1,3,16,6,21,4,17,19,15,2,23,8,13],
 [14,11,18,5,16,20,17,22,15,1,19,3,21,9,4,7,23,6,2,10,8,12],
 [13,11,12,5,9,18,7,20,10,22,3,1,6,16,4,21,8,19,2,17,23,15],
 [16,18,17,5,20,11,22,13,1,9,3,7,21,10,4,6,19,12,2,8,23,14],
 [15,18,11,5,13,20,9,22,7,1,10,3,12,6,14,4,8,21,2,19,23,17],
 [18,15,5,11,20,13,22,9,1,7,3,10,21,6,4,12,19,8,2,14,23,16],
 [17,15,16,11,13,5,12,9,14,7,10,20,6,3,8,1,4,22,2,21,23,19],
 [20,5,22,11,1,13,3,9,21,7,4,10,6,15,12,17,8,14,2,16,23,18],
 [19,5,15,11,17,13,16,9,12,7,14,10,18,6,8,3,4,1,2,22,23,21],
 [22,5,1,11,3,13,9,19,7,15,10,17,6,12,4,14,8,16,2,18,23,20],
 [21,5,19,11,15,13,17,9,16,7,12,10,14,3,6,1,8,4,18,2,20,23],
 [1,5,3,11,4,9,7,13,6,10,2,12,8,15,14,17,16,19,18,21,20,22]]
 
 
\end{BVerbatim} 
\normalsize
&
\tiny
\begin{BVerbatim}


[[3,5,4,7,6,9,8,11,10,17,13,23,15,19,12,18,16,21,14,22,20,24,2],
 [3,4,5,6,7,8,9,10,11,12,13,14,15,16,17,18,19,20,21,22,23,24,1],
 [2,1,5,23,7,21,17,24,19,22,11,15,13,18,9,16,8,12,10,14,6,20,4],
 [2,5,1,7,23,9,17,11,21,13,19,8,15,10,18,6,16,12,22,14,24,20,3],
 [2,4,1,3,23,24,21,22,17,19,7,15,11,18,13,16,9,12,8,14,10,20,6],
 [2,7,1,9,23,11,17,8,13,10,19,15,21,18,4,16,22,12,24,14,3,20,5],
 [2,6,1,4,23,3,21,24,17,22,19,5,15,18,11,16,13,20,12,14,9,10,8],
 [2,9,1,11,23,17,6,13,21,19,4,15,22,18,24,16,3,12,5,14,20,10,7],
 [2,8,1,6,23,4,17,21,11,19,13,22,15,24,18,3,16,5,12,20,14,7,10],
 [2,11,1,17,23,13,6,19,21,15,4,18,22,16,24,12,3,14,5,20,8,7,9],
 [2,10,1,8,23,6,17,4,21,9,19,24,22,3,15,5,18,7,16,20,13,14,12],
 [2,13,17,15,23,19,1,18,21,16,4,22,6,24,10,3,8,5,9,20,7,14,11],
 [2,12,17,1,23,10,6,8,21,4,19,9,22,24,15,3,18,5,16,7,20,11,14],
 [2,15,17,16,19,18,23,21,1,22,4,24,6,3,10,5,8,20,9,7,12,11,13],
 [2,14,17,12,23,1,19,6,21,10,4,8,22,9,24,13,3,11,5,7,18,20,16],
 [2,17,14,19,23,18,1,21,12,4,6,22,10,24,8,3,9,5,13,7,11,20,15],
 [2,16,14,15,12,13,1,10,23,8,6,11,4,9,21,3,24,7,22,5,19,20,18],
 [2,19,14,23,16,1,12,21,6,4,10,22,8,24,9,3,13,5,11,7,15,20,17],
 [2,18,14,16,23,12,1,15,6,10,21,8,4,13,9,11,24,3,22,7,5,17,20],
 [2,21,23,22,1,24,4,3,6,5,10,8,14,9,12,7,13,11,16,15,18,17,19],
 [2,20,23,14,1,16,12,18,6,15,10,19,8,13,4,11,9,17,3,7,24,5,22],
 [2,23,20,1,14,4,12,6,16,10,18,8,15,9,13,24,11,3,19,7,17,5,21],
 [2,22,20,21,14,18,16,19,12,15,1,13,10,17,8,11,6,9,4,7,3,5,24],
 [2,1,20,4,14,6,12,10,16,8,18,9,15,13,22,11,19,3,17,7,21,5,23]]
 
\end{BVerbatim}
\normalsize
\\
\hdashline[0.3pt/0.8pt]
\\
27-Line Solution \#1 (225 Triangles): &
27-Line Solution \#2 (225 Triangles): \\
\tiny
\begin{BVerbatim}


[[27,19,25,23,26,21,24,17,20,18,22,15,16,13,14,11,12,9,10,7,8,5,6,3,4,2],
 [3,19,5,17,7,15,13,21,11,23,9,18,16,25,14,20,12,22,6,24,8,27,10,26,4,1],
 [2,19,27,17,23,21,25,15,18,13,20,11,22,16,24,9,14,5,12,7,26,8,10,6,1,4],
 [5,19,7,17,13,15,11,21,9,23,16,18,14,25,12,20,6,22,8,24,10,27,26,2,1,3],
 [4,19,2,17,27,21,23,15,25,13,18,11,20,16,22,9,24,14,3,12,26,7,10,8,1,6],
 [7,19,13,17,11,15,9,21,16,23,14,18,12,25,20,4,22,2,24,27,8,26,10,3,1,5],
 [6,19,4,17,2,15,21,13,23,11,25,18,27,16,20,9,22,14,24,12,3,26,5,10,1,8],
 [9,13,11,19,15,17,12,16,14,21,18,23,20,25,22,4,24,2,27,6,26,3,10,5,1,7],
 [8,13,19,11,17,15,6,21,4,23,2,18,25,16,27,20,7,22,5,24,3,14,26,12,1,10],
 [11,13,12,15,14,17,16,19,18,21,20,23,22,25,24,4,27,2,26,6,3,8,5,7,1,9],
 [10,13,8,19,9,17,6,15,4,21,2,23,7,25,27,18,5,20,3,22,24,16,26,14,1,12],
 [13,10,15,19,17,8,16,21,14,23,18,6,25,4,20,2,22,27,24,7,3,5,26,9,1,11],
 [12,10,11,8,9,19,6,17,4,15,2,21,7,23,27,25,5,18,3,20,24,22,26,16,1,14],
 [15,10,17,19,16,8,21,12,23,6,18,4,25,2,20,27,22,7,24,5,3,9,26,11,1,13],
 [14,10,12,19,8,17,9,6,11,4,13,2,7,21,27,23,5,25,3,18,24,20,26,22,1,16],
 [17,10,19,14,8,12,21,6,23,4,18,2,25,9,27,7,20,5,22,3,24,11,26,13,1,15],
 [16,10,14,19,12,8,15,9,11,6,13,4,7,2,5,27,3,23,25,21,26,24,1,20,22,18],
 [19,10,21,8,23,12,6,14,4,16,2,9,25,7,27,11,5,13,3,15,24,26,20,1,22,17],
 [18,10,16,14,17,12,15,8,11,9,13,6,7,4,5,2,3,27,1,25,26,23,24,21,22,20],
 [21,10,23,8,25,6,4,12,2,14,27,9,7,16,5,11,3,13,24,15,26,18,1,17,22,19],
 [20,10,18,8,14,12,16,6,9,4,11,2,13,7,15,27,5,23,3,25,17,26,1,24,19,22],
 [23,10,25,8,4,6,2,12,27,14,7,9,5,16,3,11,24,13,26,15,1,18,17,20,19,21],
 [22,10,20,8,18,12,14,6,16,4,9,2,11,7,13,27,15,5,21,3,17,25,1,26,19,24],
 [25,10,4,8,2,6,27,12,7,14,5,9,3,16,11,22,13,20,15,18,26,17,1,21,19,23],
 [24,10,22,8,20,6,12,4,14,2,16,9,18,7,11,27,13,5,15,3,21,17,23,1,19,26],
 [27,4,2,10,6,8,3,7,5,12,9,14,11,16,13,22,15,20,18,24,17,21,1,23,19,25],
 [26,4,10,2,8,6,24,12,22,14,20,9,16,7,18,11,25,13,23,15,21,5,17,3,19,1]]
\end{BVerbatim} 
\normalsize
&
\tiny
\begin{BVerbatim}


[[2,4,3,10,6,8,5,9,7,12,11,14,13,16,15,22,18,20,17,21,19,24,23,26,25,27],
 [1,4,26,10,27,22,24,14,16,8,20,9,21,12,15,6,18,11,17,7,23,13,25,5,19,3],
 [4,1,10,26,8,22,14,27,9,16,6,20,12,24,15,21,11,18,7,17,13,23,5,25,19,2],
 [3,1,2,26,27,10,24,22,25,16,20,14,21,8,18,15,23,12,17,9,11,6,13,7,19,5],
 [6,10,8,1,9,14,12,22,11,16,7,20,15,26,18,21,13,24,17,27,23,3,25,2,19,4],
 [5,10,1,8,26,14,22,9,27,16,3,20,24,12,21,15,2,18,25,17,23,11,4,13,19,7],
 [8,10,9,1,12,14,11,22,16,5,20,26,15,27,21,24,18,3,17,2,23,25,13,4,19,6],
 [7,10,5,1,6,26,3,22,27,14,24,16,2,20,25,21,4,18,23,15,17,12,19,11,13,9],
 [10,7,1,5,14,26,22,6,27,3,16,24,20,2,21,25,15,18,12,23,17,4,11,19,13,8],
 [9,7,8,5,6,1,3,26,2,27,4,24,25,22,23,20,21,16,18,14,17,15,19,12,13,11],
 [12,1,14,7,22,5,16,26,20,27,15,24,21,3,18,2,17,25,23,6,4,9,19,8,13,10],
 [11,1,7,14,5,22,26,16,27,20,3,24,6,21,2,15,25,18,9,23,4,17,8,19,10,13],
 [14,1,16,22,15,20,18,26,21,5,24,27,17,3,23,2,25,7,4,6,19,9,8,11,10,12],
 [13,1,11,7,12,5,9,26,6,22,3,27,8,24,2,16,25,20,4,21,23,18,10,17,19,15],
 [16,1,22,13,20,5,26,7,27,11,24,3,21,6,2,12,25,9,18,4,23,8,17,10,19,14],
 [15,1,13,22,7,5,11,26,12,27,6,3,9,24,8,2,14,25,4,20,23,21,10,18,19,17],
 [18,22,20,1,21,26,24,5,27,13,3,7,2,11,25,6,23,9,4,12,8,15,10,14,19,16],
 [17,22,1,20,13,26,5,21,27,24,7,3,11,2,6,25,12,9,15,4,8,23,14,10,16,19],
 [20,22,21,1,24,26,23,27,25,3,2,5,4,7,6,13,9,11,8,12,10,15,14,17,16,18],
 [19,22,17,1,18,13,15,5,7,26,11,27,12,3,6,24,9,2,8,25,14,4,16,23,10,21],
 [22,19,1,17,26,13,5,18,27,7,24,11,3,15,6,12,2,9,25,8,4,14,23,16,10,20],
 [21,19,20,17,18,1,15,13,16,7,11,5,12,26,9,6,14,3,8,27,2,24,4,25,10,23],
 [24,1,26,19,27,5,3,13,2,7,25,11,6,17,9,12,4,15,8,18,14,21,16,20,10,22],
 [23,1,19,26,17,5,13,27,18,7,21,11,15,3,12,6,20,9,16,8,14,2,22,4,10,25],
 [26,1,27,19,3,5,2,13,7,23,11,17,6,18,12,15,9,21,8,20,14,16,4,22,10,24],
 [25,1,23,19,24,17,21,13,18,5,15,7,20,11,16,12,22,9,14,6,8,3,10,2,4,27],
 [1,25,19,23,5,17,13,24,18,21,7,15,11,20,12,16,6,9,3,14,8,22,2,10,4,26]]
\end{BVerbatim}
\normalsize
\end{tabular}
\end{adjustbox}
\end{center}
\end{samepage}

\bibliographystyle{plain}
\small
\bibliography{main}

\end{document}